\input amstex
\documentstyle{amsppt}
\magnification 1000
\catcode`\@=11
\def\logo\@{}
\catcode`\@=\active  
\NoBlackBoxes
\baselineskip=22pt
\topmatter
\title Existence and Regularity  for an Energy Maximization Problem
in Two Dimensions  \endtitle
\author Corrected Version (2015)
\endauthor
\endtopmatter
\document

\define\[{\left[}%
\define\]{\right]}%
\define\({\left(}%
\define\){\right)}%
\baselineskip=20pt

Spyridon Kamvissis

Max Planck Institute for
Mathematics in the Sciences, Leipzig, Germany

and

Department of Applied Mathematics, University of Crete, Greece 

\bigskip

Evguenii A. Rakhmanov

Department of Mathematics

University of South Florida, Tampa, Florida 33620, USA

\bigskip

Originally published in the
Journal of Mathematical Physics, v.46, n.8, 1 August 2005. 
Incorporating an addendum which appeared in the Journal of Mathematical Physics, v.50, n.9, 2009.

\bigskip

ABSTRACT

We consider  the variational problem of  maximizing the weighted
equilibrium
Green's energy of a  distribution of charges free to move in a subset of
the upper half-plane, under a particular external field.
We show that this problem admits a solution and that, under some
conditions, this solution is an S-curve (in the sense of
Gonchar-Rakhmanov).
The above problem appears in the theory of the semiclassical limit of
the integrable focusing nonlinear Schr\"odinger
equation. In particular, its solution provides a justification of a
crucial step
in the asymptotic theory of nonlinear steepest descent for the
inverse scattering problem of the associated linear non-self-adjoint
Zakharov-Shabat operator and the equivalent Riemann-Hilbert
factorization  problem.

\newpage

1. INTRODUCTION

Let 
$ \Bbb H = \{ z: Im z >0 \} $ be the complex upper-half plane  and
 $\bar \Bbb H =  \{ z: Im z \geq 0 \} \cup \{\infty\}$ 
 be the closure  of $ \Bbb H $. Let
also 
$ \Bbb K = \{ z: Im z >0 \} \setminus \{ z: Rez =0, 0< Im z \leq A \}$,
where $A$ is a positive constant.
In the closure of this space, $\bar \Bbb K $, we consider the points
$ix_+$ and $ix_-$, where $0 \leq x < A$ as distinct.
In other words, we cut a slit in the upper half-plane along the
segment $(0, iA)$ and distinguish between the two sides of the slit.
The point infinity belongs to $\bar \Bbb K$, but not $\Bbb K$.
We define $\Bbb F$ to be the set of all "continua"
$F$ in $\bar \Bbb K$
(i.e. connected compact sets)
containing  the distinguished points  $0_+, 0_-$.

Next,  let  $\rho^0(z) $ be a given complex-valued function on $ \bar \Bbb H$
satisfying 
$$ 
\aligned 
\rho^0(z)~~is~~holomorphic~~in ~~~\Bbb H,\\
\rho^0(z)~~is~~~~~~~~continuous~~~in ~~~~\bar \Bbb H , \\
Re [\rho^0(z)] =0,~~for~~~z \in [0,iA], \\
Im [\rho^0(z)] >0,~~for~~~z \in (0,iA] \cup \Bbb R.
\endaligned 
\tag1 
$$ 
Define
$G(z; \eta)$ to be the Green's function for the upper half-plane
$$
\aligned
G(z; \eta) = log {{ |z-\eta^*| } \over {|z-\eta|}}
\endaligned
\tag2
$$
and  let $d\mu^0 (\eta)$ be the nonnegative measure $-\rho^0(\eta) d\eta$
on the segment $[0,iA]$ oriented from 0 to iA. The  star denotes
complex conjugation. Let the "external field" $\phi$ be defined by
$$
\aligned
\phi (z) =
-\int G(z; \eta) d\mu^0(\eta) - Re (i\pi J \int_{z}^{iA} \rho^0 (\eta) d\eta
+2i J (z x + z^2 t) ),
\endaligned
\tag3
$$
where  $x, t$ are real  parameters with
$t \geq 0$ and $J=1,~~for~~x \geq 0,$ while
$J=-1,~~for~~x < 0$. $Re$ denotes the real part.

The particular form of this field is dictated by the particular 
application to the dynamical system
we are interested in. The conditions (1) are  natural in view of this application.
But many of our results in this paper are valid if the term
$z x + z^2 t$ is replaced by any polynomial in $z$.
Here $x, t$ are in fact the 
space and time variables for the associated PDE problem (see (9)-(10) below). 

\bigskip

Let $\Bbb M$ be  the set of all positive Borel measures on $\bar \Bbb K$, 
such  that both the free energy
$$
\aligned
E(\mu) = \int \int G(x,y) d\mu(x) d\mu(y), ~~~\mu \in \Bbb M
\endaligned
\tag 4
$$
and $\int \phi d\mu$ are finite.
Also, let
$$
\aligned
V^{\mu} (z) = \int G(z,x) d\mu(x), ~~~\mu \in \Bbb M.
\endaligned
\tag 5
$$
be the Green's potential of the measure $\mu$.

The weighted energy of the field $\phi$ is
$$
\aligned
E_{\phi} (\mu) =  E(\mu) + 2 \int \phi d\mu,  ~~~\mu \in \Bbb M.
\endaligned
\tag 6
$$

Now, given any continuum $F \in  \Bbb F$, the equilibrium measure
$\lambda^F$ supported in $F$ is defined by
$$
\aligned
E_{\phi} (\lambda^F) = min_{\mu \in  M(F)} E_{\phi} (\mu),
\endaligned
\tag 7
$$
where $M(F)$ is the set of measures in $\Bbb M$ which are supported in $F$,
provided such a measure exists. $E_{\phi} (\lambda^F)$ is the
equilibrium energy of $F$.

The aim of this paper is to prove the existence of a so-called  S-curve
([1])
joining the points $0_+$ and $0_-$ and lying entirely in $\bar \Bbb K $,
at least under some  extra assumptions.  
By S-curve we mean  an oriented  curve $F$ such that the equilibrium measure
$\lambda^F$ exists, its support consists of a finite
union of analytic arcs
and
at any interior point of $supp\mu$
$$
\aligned
{d \over {d n_+}} (\phi + V^{\lambda^F}) =
{d \over {d n_-}}  (\phi + V^{\lambda^F}),
\endaligned
\tag8
$$
where the two derivatives above denote the normal (to $supp\mu$) derivatives.

To prove the existence of the S-curve we will first need to prove
the existence of a  continuum $F$ maximizing  the equilibrium energy over 
$\Bbb F$. Then we will show that the maximizer is in fact an S-curve. 

It is not always true that  an equilibrium measure exists
for a given continuum.
The Gauss-Frostman theorem ([2], p.135) guarantees the existence
of the equilibrium measure when  $F$ does not touch the boundary
of the domain $\Bbb H$. This is not the case here.
Still, as we show in the next section,
in the  particular case of our special external field,
for any given $x,t$ and for a large class of  continua $F$ not containing
infinity, the weighted energy is bounded below and
$\lambda^{F}$ exists.
So, in particular, we do know that the
supremum of the equilibrium weighted energies over all continua is greater than
$-\infty$.

\bigskip

S-curves were first defined in [1], where the concept
first arose in  connection with the problem of
rational approximation of analytic functions.
Our own motivation comes form a seemingly completely
different problem, which is the analysis
of the so-called semiclassical asymptotics for the focusing 
nonlinear Schr\"odinger equation.  More precisely,
we are interested  in studying the behavior of solutions of
$$
\aligned
i\hbar\partial_t\psi + 
\frac{\hbar^2}{2}\partial_x^2\psi + |\psi|^2\psi = 0, \\ 
under~~~~\psi(x,0)=\psi_0(x),
\endaligned
\tag9
$$
in the so-called semiclassical limit, i.e. as $\hbar \to 0$.
For a concrete discussion, let us here assume that
$\psi_0(x) $ is a positive  "bell-shaped" function; 
in other words assume that  
$$
\aligned 
\psi_0 (x) >0, ~~~x \in \Bbb R,\\
\psi_0 (-x) = \psi_0 (x),\\
\psi_0 ~has  ~one ~single ~local ~maximum ~at ~0,~~~\psi_0 (0)=A,\\
\psi_0''(0) < 0, \\
\psi_0  ~ is ~~Schwartz.
\endaligned
\tag10
$$

This is a completely integrable partial differential equation
and can be solved via the method of
inverse scattering. The semiclassical limit is analyzed in the
recent research monograph [3].
In Chapter 8 of [3] it is noted that the semiclassical problem is related
and can be reduced to a particular "electrostatic"
variational problem of maximizing the equilibrium energy 
of a distribution of charges that are free to move 
under a  given external electrostatic field (assuming that the
WKB-approximated density of the eigenvalues admits a holomorphic extension
in the upper half-plane).
In fact, it is pointed out that the existence and regularity of
an S-curve  implies the existence of
the so-called "g-function" necessary to justify the otherwise
rigorous methods employed in [3].

We would like to point out that the problem of
the existence of the "g-function" for the semiclassical 
nonlinear Schr\"odinger problem  is 
not  a mere technicality of isolated interest. Rather, it
is an instance of a crucial element in the asymptotic
theory of Riemann-Hilbert problem factorizations associated to integrable systems.
This  asymptotic method 
has been made rigorous and systematic in [4] where in fact the
term "nonlinear steepest descent method" was first employed to 
stress the relation with the classical
"steepest descent method" initiated by Riemann in the study of
exponential integrals with a large phase parameter. 
Such exponential integrals appear in the solution of Cauchy problems for 
linear evolution equations, when  one employs the method of Fourier transforms.
In the case of nonlinear integrable equations, on the other hand, the nonlinear analog
of the Fourier transform  is the scattering transform and the inverse  problem
is now a Riemann-Hilbert factorization problem. 
While in the  "linear steepest descent method" the contour of integration must 
be deformed to a union of contours of "steepest descent" which will make the 
explicit integration of the integral  possible, in the case of the
"nonlinear steepest descent method" one deforms the original
Riemann-Hilbert factorization contour to appropriate
steepest descent contours where the resulting Riemann-Hilbert problems are
explicitly solvable. 

In the linear case, if the phase and the critical points of the phase  are real
it may not be necessary to deform the integration contour. One has rather
a Laplace integral problem on the contour given. For
Riemann-Hilbert problems 
the analog is  the self-adjointness of the underlying 
Lax operator. 
In this case  the spectrum 
of the associated linear Lax operator 
is real and the original Riemann-Hilbert contour
is real.
The  "deformation contour" 
must then stay near the real line. 
One novelty of the semiclassical problem for (9)-(10)
studied in [3] however is that, due to the non-self-adjointness
of the underlying Lax operator, the "target contour" is very specific
(if not unique) and by no means obvious. 
It is best characterized via the  solution of a
maximin energy problem, in fact it is an S-curve.
The term  "nonlinear steepest descent method" thus acquires full meaning
in the non-self-adjoint case. 

Given the importance and the recent popularity of the "steepest descent method"
and the various different applications to such topics as
soliton theory, orthogonal polynomials, solvable models in statistical
mechanics, random matrices, combinatorics and representation theory,
we believe that the present work offers an important contribution.
In particular we expect that the results of this paper may be useful in the
treatment of Riemann-Hilbert problems arising in the analysis
of general complex or normal random matrices.

On the other hand, we believe that the main results of this paper,
Theorems 3, 4, 5, 7, 8 are interesting on their own. This paper can be read without the
applications to dynamical systems in mind. It concerns existence and regularity
of a solution to an energy variational maximin problem in the complex plane.

The method used to prove the existence of the 
S-curves arising in the solution of the
"max-min" energy problem was first outlined in [1] and further developed
in [5], at least for logarithmic potentials.
But, the concrete particular problem addressed in this paper involves
additional technical issues. 

The main points of the proof of our results are:

(i) Appropriate definition of the underlying space of continua (connected compact sets)
and its topology. This ensures the compactness of our space of continua which is crucial
in proving the existence of an energy maximizing element.

(ii) Proof of the semicontinuity of the energy functional that takes
a continuum to the energy of its associated equilibrium measure (Theorem 3).

(iii) Proof of existence of an energy maximizing continuum (Theorem 4).

(iv) A discussion of how some  assumptions ensure that the maximizing continuum does not
touch the boundary of the underlying space except at a finite number of points. 
This ensures that variations of continua can be taken.

(v) Proof of  formula (22) involving the support of the equilibrium
measure on the maximizing continuum and the
external field (Theorem 5).

(vi) Proof that the support of the equilibrium
measure on the maximizing continuum consists of a union of finitely 
many analytic arcs. 

(v) Proof that the maximizing continuum is an S-curve (Theorems 7 and 8).

The paper is organized as follows.
In the rest of section 1,  
we introduce the appropriate topology for our set of
continua that will provide the necessary compactness.
In section 2,  we prove a "Gauss-Frostman" type theorem which shows that
the variational problem that we wish  to solve is not vacuous.
In section 3, we present the proof of upper semicontinuity of a particularly
defined "energy functional". In section 4, we present a proof of existence
of a solution 
of the variational problem. Existence is thus derived from the semicontinuity
and the compactness results acquired earlier. In section 5, we show that, 
at least under a simplifying assumption, the "max-min"
solution of the variational problem does not touch the boundary
of the underlying domain, except possibly at some  special points.
This enables us to eventually take variations and show that the max-min
property implies regularity of the support of the solution
and the  S-property in sections 6 and 7. 
By regularity, we mean that the support of the maximizing measure
is a finite  union of analytic arcs.
In section 8, we conclude by stating the consequence of the above results
in regard to the semiclassical limit of the nonlinear Sch\"odinger equation.

We also include three appendices. The first one discusses  in detail
some  topological facts regarding the set of closed subsets of a 
compact space, equipped with the so-called Hausdorff distance.
The fact that such a space is compact is vital for proving existence of
a solution for  the variational problem.
The second appendix presents the  semiclassical asymptotics for
the initial value problem (9)-(10) in terms of theta functions, under the 
S-curve  assumption (as in [3]). It is included so that the connection with the
original motivating problem of semiclassical NLS is made more
explicit. The third appendix shows  how to get rid of the
simplifying assumption introduced in section 5.

\newpage

Following [6] (see Appendix A.1)
we introduce an appropriate topology on $\Bbb F$.
We think of the closed upper half-plane $\bar \Bbb H$ as a compact space in the
Riemann sphere. We thus choose to equip
$\bar \Bbb H$   with the  "chordal" distance, denoted by
$d_0$,
that is  the distance between
the images of $z$ and $\zeta$ under the stereographic projection.
This induces naturally a distance in $\bar \Bbb K $ (so
$d_0 (0_+, 0_-) \neq 0).$
We also denote by $d_0$ the induced  distance between compact sets
$E, F$ in $\bar \Bbb K$:
$d_0(E,F) = max_{z \in E} min_{\zeta \in F} d_0 (z,\zeta)$.
Then, we define the so-called Hausdorff metric on the set $ I ( \bar \Bbb K   )$
of closed non-empty subsets of $ \bar \Bbb K$  as follows.
$$
\aligned
d_{\Bbb K} (A,B) = sup  ( d_0 (A,B), d_0 (B,A) ).
\endaligned
\tag11
$$

In  appendix A1, we prove the following.

LEMMA A.1. The Hausdorff metric defined by (11)   is indeed a metric.
The set $ I ( \bar \Bbb K   )$
is compact and complete.

Now, it is easy to see that  $\Bbb F$ is a closed subset of $ I ( \bar \Bbb K
)$.
Hence $\Bbb F$  is also compact and complete.

REMARKS.

1. Because of the particular symmetry $\psi(x) = \psi(-x)$ of the solution to
the Cauchy problem (9)-(10) we will restrict ourselves to the case
$x \geq 0$ from now on. We then set
$J=1$ and the external field is
$$
\aligned
\phi (z) =
-\int G(z; \eta) d\mu^0(\eta) - Re [ i\pi  \int_{z}^{iA} \rho^0 (\eta) d\eta
+2i  (z x + z^2 t) ].
\endaligned
\tag3a
$$

2. The  function $\rho^0$ expresses 
the density of eigenvalues of the Lax operator
associated to (9), in the limit as $h \to 0$.
WKB theory can be used to derive an  expression for
$\rho^0$ in terms of the initial data $\psi^0(x)$
via an Abel transform (see [3]), from which it follows that
$$
\aligned
Re [\rho^0(z)] =0,~~for~~~z \in [0,iA], \\
Im [\rho^0(z)] >0,~~for~~~z \in (0,iA].
\endaligned
$$
The rest of the conditions (1) are not a necessary
consequence of WKB theory.
In particular,
it is not a priori clear what the analyticity properties
of $\rho^0$ are.  In this paper, we $assume$, for simplicity,
that $\rho^0$  admits a continuous extension in the closed upper
complex plane  which is holomorphic  in the open upper
complex plane.
We also assume that $Im \rho^0$ is positive in the  real axis.
This will be used later to show that the maximizing continuum does not
touch the real line, except at $0_+, 0_-, \infty$. It is a
simplifying but not  essential assumption.
All conditions (1) are satisfied
in the simple case where the initial data are given by $\psi(x,0) = A sechx,$
where $A$ is a positive constant.

3. It follows that
$\phi$ is a subharmonic function in $\Bbb H $
which is actually harmonic in  $\Bbb K $ ; it  also
follows that it is upper semicontinuous in $\Bbb H $.
It is then subharmonic and  upper semicontinuous
in $\bar \Bbb H$ except at infinity.

4. Even though in the end we wish that the maximum of $E_{\phi}
(\lambda^F)$ over
"continua" $F$ is a regular curve, we will begin by studying the variational
problem over
the set of continua $\Bbb F$ and only later (in section 6) we will show that
the maximizing continuum is in fact a nice curve. The reason  is that the
set $\Bbb F$ is compact, so once we prove in section 3 the upper semicontinuity
of
the energy functional, existence of a maximizing continuum will follow
immediately.

\newpage

2. A GAUSS-FROSTMAN THEOREM

\bigskip

We claim that for any continuum
$F\in \Bbb F$, not containing the point $\infty$
and approaching $0_+, ~~ 0_-$ non-tangentially to the real line,
the weighted energy is bounded below 
and the equilibrium measure $\lambda^F$ exists. This is not true for any
external field, but it is true for the field given by (3a)
because of the
particular behavior of the function $\rho^0$ near zero.

We begin by considering  the equilibrium measure on the 
particular contour $F_0$ that wraps 
itself around the straight line segment $[0,iA]$, say $\lambda^F_0$. We have

PROPOSITION 1. Consider the contour $F_0 \in \Bbb F$ consisting of the straight line segments
joining $0_+$ to $iA_+=iA$ and $iA=iA_-$ to $0_-$. The equilibrium measure
$\lambda^F_0$ exists. Its support is the imaginary segment
$[0, ib_0(x)]$, for some $0 < b_0(x) \leq A$, lying on the right of the slit
$[0, iA]$.  It can be written as $\rho (z) dz$ where $\rho (z)$  is
a differentiable   function in $[0, ib(x)]$.

PROOF: See section 6.2.1 of [3]; $\rho (z)$ can be expressed explicitly when t=0. But
note that the field $\phi$ is independent of
time on $F_0$, so $\lambda^F_0$ is also  independent of time.

\bigskip

From Proposition 1,
it follows that the maximum equilibrium energy over continua is bounded below.
$$
\aligned
max_{F \in \Bbb F} E_{\phi} (\lambda^F) = 
max_{F \in \Bbb F} min_{\mu \in  M(F)} E_{\phi} (\mu) > -\infty.
\endaligned
\tag12
$$

The following formula is easy to verify.
$$
\aligned
E_{\phi} (\mu) - E_{\phi} (\lambda^F) = E(\mu-\lambda^F) 
+2 \int (V^{\lambda^F} + \phi ) d(\mu-\lambda^F),
\endaligned
\tag13
$$
for any  $\mu$ which is  a positive measure on  the continuuum $F$.
Here 
$$
\aligned
V^{\lambda^F} (u)= \int G(u,v) d\lambda^F(v),
\endaligned
$$
where again $G(u,v)$ is the Green function for the upper
half-plane.

To show that $E_{\phi} (\mu) $ is bounded below, all we need to show
is  that the difference
$E_{\phi}(\mu)-E_{\phi} (\lambda^F)$ is bounded below.

Note that since $V^{\lambda^F} + \phi =0,~~~on~~~supp(\lambda^F)$, the integral
in 
(13) can be written
as $\int (V^{\lambda^F} + \phi) d\mu$.

We have
$$
\aligned
V^{\lambda^F} (z)+ \phi (z)= 
\int^{b_0(x)}_0 log {{|z+iu|}\over{|z-iu|}} (- \rho)_{t=0} du + \phi= \\
= - Re [ \int_0^{b_0(x)} log {{|z+iu| } \over {|z-iu|}} u^{1/2} du] + 
O(|z|) = \\
=O(|z|)~~~~~near~~~z=0. 
\endaligned
$$
So we can  
write $V^{\lambda^F} + \phi \geq  c(A,x)  |z|$ 
in a neighborhood of $z=0$,
where $c(A,x)$ will be some negative constant independent of $z$.
Note that the dependence on $t$ is not suppressed, but it
is of order $O(|z^2|).$

It is now not hard to see that the $O(|z|)$ decay implies our result, 
at least if we suppose that $F$ is contained in some sector $ \pi < \alpha < arg(\lambda) < \beta <0$
as $\lambda \to 0$. 

Write $\mu = M \sigma$, where $M>0$ is the total mass of $\mu$ and $\sigma$ is a
probability 
measure (on $F$). 
Choose $\epsilon$ such that for $|u| < \epsilon$ we have
$V^{\lambda^F} + \phi \geq  c(A,x) |u|$.
Then

$$
\aligned
E_{\phi} (\mu) - E_{\phi} (\lambda^F) \geq 
\int G(u,v) d(\mu-\lambda^F)(u)d(\mu-\lambda^F)(v) \\
+ 2 \int_{|v| \geq \epsilon} (V^{\lambda^F} + \phi)(v) d\mu(v)
+\int_{|v| < \epsilon} 2 c(A,x) |v| d(\mu-\lambda^F) (v).
\endaligned
\tag14
$$

The first integral  of the right hand side (RHS) can be written
as $\int_{|u| \geq \epsilon, |v| \geq \epsilon} + 2 \int_{|u| < \epsilon, |v|
\geq \epsilon}
+ \int_{|u|,|v| <\epsilon}. $
The sum of the first integral plus the second term of the RHS of (14) is bounded
 below, by 
the standard Gauss-Frostman theorem ([6], p.135).
It remains to consider

$$
\aligned 
(\int_F +\int_{|v| \geq \epsilon} )[ \int_{|u| < \epsilon} G(u,v)
d(\mu-\lambda^F) (u) ] 
d(\mu-\lambda^F) (v) +  \int_{|v| < \epsilon} 2 c(A,x) |v|  d(\mu-\lambda^F)(v)
\\
\geq (\int_F + \int_{|v| \geq \epsilon})
[\int_{|u| < \epsilon} G(u,v) d(M \sigma-\lambda^F) (u) +2 c(A,x)  |v| ] 
d(M \sigma-\lambda^F)(v).\\
\endaligned
\tag15
$$
Now,  it is easy to see that  since $F$ is non-tangential to 
the real line, $G(u,v) \geq const.~max\{sin(\alpha), sin(\beta)\}$ and so 
for $M=M(\epsilon)$ large enough (e.g. $M=O(\epsilon^{-2}$)
$$
\aligned
\int_{|u| < \epsilon} G(u,v) d(M\sigma - \lambda^F) (u) +2 c(A,x) |v| \geq 
large~positive~constant  + 2c(A,x)  |v| .
\endaligned
\tag16
$$
Hence the integral in (16) is positive.
Integrating again with respect to
$M d\sigma -d\lambda^F$, again for $M$ large, we see that the integral of (15)
is positive.

Since  for $M$  bounded above we have
our estimates trivially, we clearly get 
boundedness below over the set of all positive $M$.

We have thus proved one part of our (generalised) Gauss-Frostman Theorem.

THEOREM 1. Let $\phi$ be given by (3a).
Let $F$ be a continuum in $\bar \Bbb K \setminus \infty$ and
suppose that $F$ is contained in some sector $ \pi < \alpha < arg(\lambda) <
\beta <0$
as $\lambda \to 0$. Let
$M(F)$ be the set of  measures $\in \Bbb M$ which are supported in $F$.
(So, in particular their free energy is finite and 
$\phi \in L_1 (\mu)$.)
We have 
$$
\aligned
inf_{\mu \in M(F)} E_{\phi} (\mu) > -\infty.
\endaligned
\tag17
$$
Furthermore  the equilibrium measure on $F$ exists, that is there is a
measure $\lambda^F \in M(F)$ such that 
$E_{\phi} [F] = E_{\phi} (\lambda^F)= inf_{\mu \in M(F)} E_{\phi} (\mu).$ 

PROOF: The proof that (17) implies the existence of
an equilibrium measure is  a well known theorem. For our particular field
$\phi$ given by (3) it is easy to prove. Indeed, the identity
$$
\aligned
E(\mu-\nu) = 2 E_{\phi} (\mu) + 2 E_{\phi} (\nu) - 4 E_{\phi} ({{\mu+\nu} \over
2})
\endaligned
$$
implies that  any sequence $\mu_n$ minimizing $ E_{\phi} (\mu)$ is a Cauchy
sequence
in (unweighted) energy. Since the space of positive measures  is complete 
(see for example [7], Theorem 1.18, p.90),
there is a measure $\mu_0$ such that $E(\mu_n-\mu_0) \to 0$. We then have 
$E(\mu_n) \to E(\mu_0) < +\infty$ and hence $\mu_n \to \mu_0$ weakly 
(see e.g. [7], p.82-88; this is a standard result).

The fact that $\phi \in L_1 (\mu_0)$ is trivial for our particular field.

\newpage

3. SEMICONTINUITY OF THE ENERGY FUNCTIONAL

\bigskip

Let $F$ be a continuum 
contained in some sector $ \pi < \alpha < arg(\lambda) <
\beta <0$
as $\lambda \to 0$. 
We consider the functional that takes  $F$ to its  equilibrium
energy:
$$
\aligned
\Bbb E: F \to E_{\psi} [F] = E_{\psi} (\lambda^F) = \inf_{\mu \in M(F)} (E(\mu)
+ 2 \int \psi d\mu)
\endaligned
\tag 18
$$
and we want to show that it is continuous, if $\psi$ is continuous in $\bar \Bbb
H $.
Note that this is not the case for the field $\phi$ given by (3a), since it has a
singularity at $\infty$; that field
is only upper semicontinuous.
We will see how to circumvent this difficulty later. For the moment, 
$\psi$ is simply assumed to be a continuous function in $\bar \Bbb H$.

\bigskip

THEOREM 2. If $\psi$ is a continuous function in $\bar \Bbb
H \setminus \infty$ then the energy functional defined by (18) is continuous 
at any given  continuum $F$ contained in the sector $ \pi < \alpha < arg(\lambda) <
\beta <0$
as $\lambda \to 0$, not containing the point $\infty$.

PROOF: Suppose $G \in \Bbb F$, 
with $d_{\Bbb K} (F, G) < d$, a  small
positive constant such that $G$ is also contained in the sector $ \pi < \alpha < arg(\lambda) <
\beta <0$
as $\lambda \to 0$, not containing the point $\infty$.

Let $\lambda = \lambda^F_{\psi} $ be the equilibrium measure on $F$ and
$\mu = \lambda_{\psi}^G$ be the equilibrium measure on $G$. 

We consider the Green's balayage of $\mu$ on $F$, say $\hat \mu$.
Then $supp ~\hat \mu \in F$ and

$$
\aligned
V^{\hat \mu}= V^{\mu} ~~~on~~~~F,\\
\int u d\hat \mu = \int u d\mu,
\endaligned
$$
for any function $u$ that is harmonic in $\Bbb H \setminus F$ and continuous in
$\bar \Bbb H$.  

Similarly consider $\hat \lambda $, 
the balayage of $\lambda $ to $G$. We trivially have
$$
\aligned
E_{\psi} [G] \leq E_{\psi} (\hat \lambda),\\
E_{\psi} [F] \leq E_{\psi} (\hat \mu).
\endaligned
$$

LEMMA 1. Suppose $Q \in \Bbb F$,
$\mu$ some positive measure supported in $\Bbb K$ 
and $\hat \mu$ is the Green's balayage to $Q$.
Then
$$
\aligned
V^{\hat \mu} =  V^{\mu}- V^{\mu}_{Q^c},\\
E(\hat \mu) = E(\mu) - E_{Q^c} (\mu),
\endaligned
$$
where  $E_{Q^c} (\mu) $ is the unweighted Green energy with respect 
to  $Q^c = \Bbb K  \setminus Q$.
In particular, since unweighted energies are nonnegative,
$$
\aligned
E(\hat \mu) \leq E(\mu).
\endaligned
$$

PROOF: The first identity follows from the fact that
$V^{\hat \mu} - V^{\mu}  $ vanishes on $Q$ and the real line,
and is harmonic in $Q^c$ and superharmonic
in $\Bbb K$.

Integrating
$E(\hat \mu) = \int V^{\hat \mu} d\hat \mu = 
\int V^{\mu} d\hat \mu - \int V^{\mu}_{Q^c} d\hat \mu =
\int (V^{\mu} - V^{\mu}_{Q^c} ) d\mu =
E(\mu) - E_{Q^c} (\mu).$ The proof of the Lemma follows.

\bigskip

So, let $u_{\psi}$ be a function harmonic in $\Bbb H \setminus
F$ 
such that $u_{\psi} = \psi$ on $F$ and $u_{\psi} =0$ on $\partial \Bbb H
\setminus F$.
By the definition of balayage one has $\int \psi d\hat \mu = \int u_{\psi}
d\mu$.

We have $E_{\psi} [F] \leq E_{\psi} (\hat \mu) =
E(\hat \mu)  + 2 \int \psi d\hat \mu  \leq  
E(\mu) + 2 \int \psi d\mu + 2 \int (u_{\psi} - \psi) d\mu =
E_{\psi} (\mu) + 2 \int (u_{\psi}- \psi ) d\mu =
E_\psi [G] + 2 \int (u_{\psi} -\psi ) d\mu.$

In a small neighbourhood of $F,  
\bar F_{d} = \{ z:d(z,F) \leq d \}$, we have
$$
\aligned
| 2 \int (u_{\psi} - \psi) (y) d\mu(y) | \leq
C max_{y \in \bar F_{d}} |u_{\psi} (y) - \psi (y) |.
\endaligned
\tag19
$$

We assumed here that the equilibrium measures on continua near 
$F$ are bounded above. This is easy to see. Suppose, first, that
the point $\infty$ is not in $F$. Indeed,
on  the support of the equilibrium measure $\lambda$, we have
$$
\aligned
V^{\lambda} + \psi = 0.
\endaligned
$$
If the equilibrium measures on continua near $F$ were unbounded, then so
would be  the potentials $V^{\lambda} $. 
(This follows easily from explicit formulae for the
equilibrium measures in terms of the potentials.) But $\psi$ is
definitely bounded near $F$. This contradicts the above equality.

\bigskip
Now given $y \in \bar F_{d}$, choose $z \in F$ such that $|z-y| = d$.
The above expression (19) is less or equal than
$$
\aligned
 C max_{y \in \bar F_{d}} |u_{\psi} (y) - \psi (y) 
-u_{\psi} (z) + \psi (z) |
\leq o(1) + C max_{y \in \bar F_{d}} |u_{\psi} (y) - u_{\psi} (z)|.
\endaligned
$$
It remains to bound $|u_{\psi} (y) - u_{\psi} (z)|$ by an $ o(1)$ quantity.

\bigskip

The next Lemma  is due to Milloux and can be found in [8].

LEMMA 2. Suppose $D$ is an open disc of radius $R$, 
with center $z_0$; let $y$
be a  point in $D$, $F$  a continuum in $\Bbb C$, containing
$z_0 $, and $\Omega$ be the
connected component of $D \setminus F$ containing 
$y$. Let $w(z)$ be a function harmonic 
in $\Omega$ such that
$$
\aligned
w (z) =0, ~~~z \in F \cap \partial \Omega,\\
w( z) =1, ~~~~z \in \partial \Omega \setminus F.
\endaligned
$$
Then $w(y) \leq C ({{|y-z_0|} \over R})^{1/2}$.

PROOF: See [8], p.347.

\bigskip

Now, select a disc of radius $d^{1/2}$, centered on $z$.
We have 
$ |u_{\psi} (y) - u_{\psi} (z) | = o(1) $ on the part of $F$ lying in the disc,
while $|u_{\psi} (y) - u_{\psi} (z) | $ is bounded by some positive constant $M$
on
the disc boundary.

LEMMA 3. Let $\Omega$ be a domain, $\partial \Omega= F_1 \cup F_2$ and
$$
\aligned
w_1 =0,~~~ z \in F_1,\\
=1, ~~~ z \in F_2;\\
w_2 = 1, ~~~ z \in F_1, \\
=0, ~~~ z \in F_2.
\endaligned
$$
Suppose $u$ is harmonic in $\Omega$ and
$$
\aligned
u(z) \leq \epsilon, ~~~ z \in F_1, \\
u(z) \leq M, ~~~~~~~~~~~ z \in F_2.
\endaligned
$$
Then $u(z) \leq \epsilon w_2 (z) + M w_1 (z)$.

PROOF: Maximum principle.

\bigskip

Now, using Milloux's Lemma, we get
$|u_{\psi} (y) - u_{\psi} (z)| \leq o(1) w_2(z) + M w_1(z) \leq
o(1) + M C {{ |y-z| } \over d}^{1/2} \leq o(1) + M C d^{1/2}.$
This concludes the proof of Theorem 2.

\newpage

We now recall that the energy continuity proof was based on the continuity of
$\psi$.
In our case, $\phi$ is upper semicontinuous and discontinuous at $\infty$. Still
we can prove that the energy
is upper semicontinuous and that will be enough.

THEOREM 3. For the external field given by (3a), 
the energy functional defined in
(18) is upper semicontinuous on $\Bbb F^{\beta}_{\alpha}$ which consists
of continua $F$  contained in the sector $ \pi < \alpha < arg(\lambda) <
\beta <0$
as $\lambda \to 0$.

PROOF: We first note that if the external field $\phi' $ is
upper semicontinuous away from  infinity 
then so is the energy functional that takes
a given continuum $F$ to the equilibrium energy of $F$.
Indeed, if $\phi'$ is
upper semicontinuous away from  infinity, then there exists a sequence  of 
continuous functions (away from  infinity)
such that $\phi_n \downarrow \phi'$.  Each functional $E_{\phi_n} [F]$ is
continuous, away from infinity,
and $E_{\phi_n} [F] \downarrow E_{\phi'} [F]$. So, $E_{\phi'} [F]$ is upper
semicontinuous, away from infinity.

Now consider the field $\phi$ given by (3a).
Let $F$ be a continuum. If $\infty$ is not in $F$, then we're done.
If $\infty \in F$, let $\lambda=\lambda^F$ be the equilibrium measure.
We can assume that on the equilibrium measure 
$\phi$ is bounded by 0.
Indeed,
on  the support of the equilibrium measure $\lambda$, we have
$$
\aligned
V^{\lambda} + \phi = 0.
\endaligned
$$
But $V^{\lambda} \geq 0$, so $\phi \leq 0$.

This means that we can change $\phi$ to $\phi'=min(\phi, 0),$
which $is$ an upper semicontinuous function. Theorem 3 is proved.

\bigskip

REMARK. If we naively consider the  functional taking a measure to its
weighted energy we will see that it is not continuous even if the external field
is continuous. It is essential that the energy functional is defined on
equilibrium measures.

\newpage

4. PROOF OF EXISTENCE OF A MAXIMIZING CONTINUUM
 
\bigskip 

THEOREM 4.  For the external field given by (3a),
there exists a continuum  $F \in \Bbb F^{\beta}_{\alpha}$ such that the equilibrium measure
$\lambda^F$ exists and 
$$
\aligned 
E_{\phi} [F] (= E_{\phi} (\lambda^F)) =
max_{F \in \Bbb F} min_{\mu \in M(F)} E_{\phi} (\mu). 
\endaligned
$$
PROOF: We know (see for example section 2)
that there is at least one continuum $F$
for which the equilibrium 
measure exists and $E_{\phi} (\lambda^F) > -\infty$, for all time. 
On the other hand, clearly $E_{\phi} (\lambda^F) \leq 0$ for any $F$. Hence the
supremum over continua in $\Bbb F$ is finite (and trivially nonpositive), 
since $\Bbb F$ is compact. Call it $L$.
 
We can now take a  sequence $F_n$ such that $E_{\phi} [F_n] \to L$. 
Choose  a convergent
subsequence of continua $F_n \to F$, say. 
By upper semicontinuity of the weighted energy functional,
$$
\aligned
limsup E_{\phi} [F_n]
\leq E_{\phi} [F] \leq L = lim E_{\phi} [F_n].
\endaligned
$$
So $L= E_{\phi} [F]$. 
The theorem is proved.

\bigskip

5. ACCEPTABILITY OF THE CONTINUUM

\bigskip

We have thus shown that a solution of the maximum-minimum problem exists.
We do not know yet that the maximizing continuum  is a contour. 
Clearly the pieces of the continuum lying in the region where the external field is positive 
do not support the equilibrium measure and by the continuity of the external field they can be 
perturbed to a finite union of analytic arcs. The real problem is to show that
the support of the equilibrium measure is a finite  union of analytic arcs. This will follow
from the analyticity properties of the external field.

Note that the maximizing continuum cannot be unique, since the subset
where the equlibrium measure is zero can be perturbed without changing the
energy. A more interesting question is whether the support of the equilibrium measure 
of the maximizing contour is unique. We do not know the answer to this question
but it is not important as far as the application to the semiclassical limit
of the nonlinear Schr\"odinger equation is concerned. (See Appendix A2.)

It is important however, that the maximizing continuum does not approach  the
boundary of the underlying space  except of course at the points $0^-, 0^+$,
and perhaps at $\infty$.
This is to guarantee that variations with respect to the maximizing
contour can be properly taken.

The proof of the acceptability of the continuum requires two things.

(i) The continuum does not approach the real negative axis.

(ii) The continuum does not approach the real positive axis.

We will also  make the folowing assumption.

ASSUMPTION (A). The continuum maximizing the equilibrium energy
does not touch the linear segment $(0, iA]$.
 
REMARK. Assumption (A) is not satisfied at $t=0$, where in fact the 
continuum is a contour $F_0$ wrapping around the linear segment $[0, iA]$. 
However, the case $t=0$ is well understood. The equilibrium measure for
$F_0$ exists and its support is  connected.
On the other hand
assumption (A) $is$ satisfied for small $t >0$. (See Chapter 6 of [3].)

REMARK. It is conceivable that at some positive $t_0$ there is an $x$ for which 
assumption (A) is not satisfied. It can in fact be dropped but the analysis of the semiclassical
limit of NLS will get more tedious; see Appendix A3.

\bigskip

PROPOSITION 2.
The continuum maximizing the equilibrium energy
does not approach the real axis except at the points zero and possibly infinity.
More precisely, if the  real positive numbers $\alpha <\beta <\pi$ are small,
then it 
does not touch the boundary of $\Bbb F^{\beta}_{\alpha}$ near $0$ nor the real axis.

PROOF: (i) If $z <0$, then $\phi (z) = \pi \int_z^0 |\rho^0 (\eta)| d\eta >0$. 

This follows from an easy
calculation, using the conditions defining $\rho^0$.
But we can always delete 
the strictly positive measure lying in a region where the field 
is positive and make the energy smaller. So even  the solution of the "inner"
minimizing problem must lie away from the real negative axis.

(ii) If $z >0$, then again
a short calculation shows that ${{d\phi }\over {dImz}}  >0,~~~for~~~t>0 $.

It is crucial here that if $u \in \Bbb R$ then
$G(u,v)=0$, while if both $u, v $ are off the real line 
$G(u,v) > 0.$  Hence, for any configuration that involves a
continuum  including points on the real line, we can find a configuration
with no points on the real line, by pushing
measures up away from the real axis, which  has
greater  (unweighted $and$ weighted) energy.  So,
suppose the maximizing continuum touches the axis. We can always push the
measures up away from the real axis and end up with a continuum that has
greater minimal energy, thus arriving at a contradiction.

The proposition is now proved.

\bigskip

REMARK.
It is also important  to consider the point at infinity. We cannot prove that
the continuum does not hit this point. (In fact, our numerics
([3], Chapter 6) show that it may well do so.)
In connection with the semiclassical problem (9)-(10)
as analyzed in [3], it might seem at first that the maximizing continuum 
should not pass through infinity.
Indeed, the transformations (2.17) and (4.1) of [3] implicitly assume that
the continuum $C$ lies in $\Bbb C$. Otherwise, one would lose the
appropriate normalization for $M$ at infinity. However, one must 
simply notice that infinity is just an arbitrary choice of normalizing
point, once we view our Riemann-Hilbert problems in the compact Riemann
Sphere. The important observation is that the 
composition of transformations (2.17) and (4.1) 
(which are purely formal, i.e. no estimates are required and no approximation
is needed) does not introduce any bad (essential) singularities. In the end,
the asymptotic behavior of $\tilde N^{\sigma}$ is still the identity as 
$z \to \infty$ in the lower half-plane  and non-singular as $z \to \infty$
in the upper half-plane.
So, in the end it $is$ 
acceptable for a continuum  to go through the point infinity.

\newpage

6. TAKING SMALL VARIATIONS

\bigskip

We now  complexify the external field and extend it to
a  function in 
the whole complex plane, by turning a Green's potential to a logarithmic potential.
We will thus be able to make direct use of the
results of [5].

We let, for any complex $z$,
$$
\aligned
V(z)= -\int_{-iA}^{iA} log(z-\eta) \rho^0(\eta) d\eta
- (2ixz+2itz^2 +i\pi \int_{z}^{iA}  \rho^0(\eta) d\eta)
\endaligned
\tag20
$$
and $V_R = Re V$ be the real part of $V$. 
In the lower half-plane the function $\rho^0$ is
extended simply by 
$$
\aligned
\rho^0 (\eta^*) = (\rho^0 (\eta))^*.
\endaligned
\tag20a
$$
Note right away that the field $\phi $ defined in (3a) is the
restriction of $V_R = Re V$ to the closed upper half-plane.

The actual contour of the logarithmic integral is chosen to be the linear
segment
$\chi$ joining the points $-iA, 0, iA$. 
The branch of the logarithm function $ log(z-\eta) $ is then defined  to agree
with the principal 
branch as $z \to \infty$ , and with jump across  the very contour $\chi$.

The unweighted Green's energy (4) can be written as
$$
\aligned
E_{V_R} (\mu) =  \int_{supp\mu} \int_{supp\mu} log {1 \over {|u-v|}}
d\mu(u) d\mu(v) +  2 \int_{supp\mu} V_R (u) d\mu(u),
\endaligned
\tag21
$$
where the measures $\mu$ are extended to the lower half complex plane by 
$$
\aligned
\mu(z^*) = -\mu(z).
\endaligned
\tag21a
$$
(So they are "signed" measures.)

Having established in section 5 that the contour solving the variational
problem does not touch the boundary of the underlying space except at three  specific points,
we can take
small variations  of measures and contours, never intersecting that boundary, and keeping the points $0_-, 0_+$ fixed. 
In view of (21a) we can think of them as variations of measures symmetric under
(21a) in the full complex plane, 
never approaching  the
real line, and keeping the points $0_-, 0_+$ fixed.
The perturbed  measures  do not
change sign. The fact that $\infty$
can belong to the contour is not a problem. Our variations will keep it
automatically fixed.

The first step is to show that the solution of the variational problem 
satisfies  a crucial relation.

\bigskip

REMARK. It is not hard to see that
the variational problem of Theorem 4 is actually $equivalent$ to the variational 
problem of maximizing equilibrium measures on continua in the whole complex plane,
under the symmetry (21a) and the condition that measures are to  positive 
in the upper half-plane and negative in the lower half-plane.

\bigskip

THEOREM 5. Let $F$ be the maximizing continuum of Theorem 4
and $\lambda^F$ be the  equilibrium measure minimizing the
weighted logarithmic energy (6)  under the
external field $V_R = ReV$ where $V$ is
given by (20). Let $\mu$ be the extension of $\lambda^F$ to the
lower complex plane via $\mu(z^*) = -\mu(z)$.
Then
$$
\aligned
\int {{d\mu(u)} \over {u-z}} + V'(z) )^2= 
 V'(z))^2 - 2  \int {{V'(z)-V'(u)} \over {z-u}} d\mu(u)  \\ 
+   {1 \over z^2}  \int 2 (u+z) V'(u) ~ d\mu(u)   .
\endaligned
\tag22
$$

PROOF: We first need to prove the following.

THEOREM 6. Let $\Gamma$  be a critical point of
the functional taking a continuum $\Gamma \in \Bbb F $ to 
$E_{V_R} (\lambda^{\Gamma})$, and assume that 
$\Gamma$ is not  tangent to
$\Bbb R $. Also assume that $\Gamma$ does not touch the segment
$[0, iA]$ except at zero.  Let $\mu$ be the extension of $ \lambda^{\Gamma}$
via $\mu(z^*) = -\mu(z)$,
 $O_{\Gamma}$ be an open set containing the interior of 
$\Gamma  \cup \Gamma^* $ and
$h \in C^1(O_{\Gamma} )$ such that $h(0) =0$. We have
$$
\aligned
 \int \int   {{ h(u)-h(v)} \over {u-v}}  d\mu(u) d\mu(v) = 2 \int V'(u) h(u) d\mu(u) .
\endaligned 
\tag23
$$
PROOF: Consider the family of (signed) measures
$\{\mu^{\tau}, \tau \in \Bbb C, |\tau| < \tau_0 \}$
defined by $d\mu^{\tau} (z^{\tau}) =d\mu(z)$ where
$z^{\tau} = z + \tau h(z))$, or equivalently,
$\int f(z) d\mu^{\tau} (z) = \int f(z^{\tau}) d\mu (z),~~~
f\in L_1 (O_{\Gamma})$.
Assume that
$\tau$ is small enough (so that the support of the deformed continuum does not
hit the linear segment $(0, iA]$ and does not come close to the real line near $0$ except at $0$).

With $\hat h = \hat h (u,v) =  {{ h(u)-h(v)} \over {u-v}}$,
we have $ {{ u^{\tau}-v^{\tau}} \over {u-v}} = 1 +\tau \hat h$,
so that $log {1 \over { |u^{\tau}-v^{\tau}|} } -
log {1 \over { |u-v|} } = - log|1 +\tau \hat h| = -Re (\tau \hat h)
+ O( \tau^2).$

Integrating with respect to $ d\mu(u) d\mu(v)$ we arrive at
$$
\aligned
E(\mu^{\tau}) - E(\mu) = -Re [\tau \int \hat h  d\mu(u) d\mu(v)]
+ O( \tau^2),
\endaligned
\tag24
$$
where $E(\mu)$ denotes the free logarithmic energy of the measure $\mu$. 
Also, 
$$
\aligned
\int V_R d\mu^{\tau} - \int V_R d\mu =
2 \int (V_R(u^{\tau}) - V_R(u)) d \mu (u) =\\
2Re [ \tau \int V'(u) h(u) d\mu(u) ] 
+ O (\tau^2).
\endaligned
$$
Combining with the above,
$$
\aligned
E_{V_R}(\mu^{\tau}) - E_{V_R}(\mu) = 
Re (-\tau \int \hat h  d\mu(u) d\mu(v) + 2 \tau \int V' h d\mu )
+ O( \tau^2).
\endaligned
\tag25
$$

So, if $\mu$ is (the symmetric extension of)  a critical point 
of the map $\mu \to E_{V_R} (\mu)$ the linear part of
the increment is zero. 
In other words  given a $C^1$ function $h$
and a measure $\mu$ the function $E_{V_R}(\mu^{\tau})$ 
of $ \tau$ is differentiable at $\tau=0$ and the derivative is
$$
\aligned
Re (-H(\mu)),~~~where~~
H(\mu)= \int\int \hat h d\mu^2 - 2\int V' h d\mu.
\endaligned
\tag26
$$

But what we really want is the derivative of the energy as a function of
the equilibrium measure. 
This function can be shown to be differentiable and  its derivative 
can be set to zero at a critical continuum.

Indeed, we need to show the following.

LEMMA 4. 
$$
\aligned
{d \over {d\tau}} E_{V_R} ((\lambda^{\Gamma})^{\tau}) |_{\tau=0} =
{d \over {d\tau}} E_{V_R} (\lambda^{{\Gamma}_{\tau}} ) |_{\tau=0} =0.
\endaligned
$$
In the relation above
$ \Gamma_{\tau} = supp (\lambda^{\Gamma})^{\tau}$.
The first derivative is of a function of general measures. The second
derivative is of a function of equilibrium measures.

PROOF: Define the measure $\sigma_{\tau}$ with support $\Gamma$ and such
that $(\sigma_{\tau})^{\tau} = \lambda^{\Gamma_{\tau}}$.  

LEMMA 5.  With $H$ defined by (26), we have
$$
\aligned
Re H(\sigma_{\tau}) \to Re H(\lambda^{\Gamma}),
\endaligned
$$
as $\tau \to 0.$

PROOF. By (25)-(26), we have
$$
\aligned
E_{V_R} ((\lambda^{\Gamma})^{\tau}) - E_{V_R} (\lambda^{\Gamma}) = 
-Re (\tau H(\lambda^{\Gamma}) +O(\tau^2)),\\
E_{V_R} (\lambda^{\Gamma_{\tau}}) - E_{V_R} (\sigma_{\tau}) =
-Re (\tau H(\sigma_{\tau}) +O({\tau}^2)).
\endaligned
$$
On the other hand,
$E_{V_R} (\sigma_{\tau}) \geq E_{V_R} (\lambda^{\Gamma}), $ and
$E_{V_R} ((\lambda^{\Gamma})^{\tau}) \geq E_{V_R} (\lambda^{\Gamma_{\tau}}).$
It follows that
$$
\aligned
E_{V_R}(\sigma_{\tau}) - Re (\tau H(\sigma_{\tau})) + 
O(\tau^2) =E_{V_R}(\lambda^{\Gamma_{\tau}}) \leq
E_{V_R}((\lambda^{\Gamma})^{\tau}) = E_{V_R}(\lambda^{\Gamma}) -
Re (\tau H(\lambda^{\Gamma})) +O(\tau^2).
\endaligned
$$ 
Hence $E_{V_R}(\sigma_{\tau}) \to E_{V_R}(\lambda^{\Gamma})$.  

As in the proof of Theorem 1,
it follows that $\sigma_{\tau} \to \lambda^{\Gamma}$ weakly;
see [7], pp.82-88.
It  then  follows  immediately 
that  $ H(\sigma_{\tau} ) \to  H(\lambda^{\Gamma}).$ 
This proves Lemma 5.

To complete the proof of Lemma 4, we note that
$0 \geq E_{V_R}((\lambda^{\Gamma})^{\tau})-E_{V_R}(\lambda^{\Gamma}) 
\geq E_{V_R}(\lambda^{\Gamma_{\tau}}) -E_{V_R}(\sigma_{\tau}) = 
-Re(\tau H(\sigma_{\tau})) + O(\tau^2)$. 
Hence the derivative of $E_{V_R} ((\lambda^{\Gamma})^{\tau})$
at $\tau =0$ is equal to the derivative of $E_{V_R} (\lambda^{\Gamma_\tau})$
at $\tau =0$ which is equal to 
$Re H(\lambda^{\Gamma}) $. This proves Lemma 4 and Theorem 6,
by considering both $\tau$ real and $\tau$ imaginary.

PROOF OF THEOREM 5. Consider the Schiffer variation, i.e. take
$h(u) = {{ u^2} \over {u-z}}$ where $z$ is some fixed point not in $\Gamma$.
Note that $h(0) =0$ so that the deformation
$z^{\tau} = z + \tau h(z)$ keeps the points $0_+, 0_-$ fixed. Also assume that
$\tau$ is small enough so that the support of the deformed continuum does not
hit the linear segment $(0, iA]$ or a non-zero point in the real line.
We have
$$
\aligned
\hat h = \hat h (u,v) =  {{ h(u)-h(v)} \over {u-v}}=
1-{{ z^2 } \over {(u-z)(v-z)}},
\endaligned
$$
and therefore
$$
\aligned
 \int \int \hat h (u,v) d\mu(u) d\mu(v)  =  
\int \int  d\mu(u) d\mu(v)
- z^2 [ \int_{supp\mu} {{d\mu(u)} \over {u-z}} ]^2 .
\endaligned
$$

Next, we have
$$
\aligned
 \int 2 V'(u) h (u) d\mu(u)  \\ = 2 \int (u+z)   V'(u) d\mu(u) 
+  z^2  \int  { V'(u){d\mu(u)} \over {u-z}}  \\ = 
 \int 2 V'(u)  (u+z) d\mu(u)  + \\
 2 z^2  \int  { {V'(u) - V'(z) } \over {u-z}} d\mu(u) 
+2 z^2 V'(z)   \int {{d\mu(u)} \over {u-z}}.  
\endaligned
$$
Theorem 5 now follows from Theorem 6.

\bigskip

REMARK. If our  continuum is allowed to touch the point $ iA$
(so we slightly weaken assumption (A)) then we may
need to  keep  points $\pm iA$ fixed under a small variation. We can then choose the Schiffer variation
$h(u) = {{ u^2 (u^2 +A^2)} \over {u-z}}$. We will arrive at a similar and equally useful formula.

In general if one  wants to keep points $a_1, ... , a_s$ fixed, the appropriate Schiffer variation is
$h(u) = {{ \Pi_{i=1}^{i=s} (u-a_i)} \over {u-z}}$.

\bigskip

PROPOSITION 3. The support of the equilibrium measure
consists of a finite number of analytic arcs.

PROOF:
Theorem 5 above implies  that the support of $\mu$ is the level set of 
the real part of a function that is analytic except at 
countably many branch points. 
In fact,  $supp \mu$ is characterized by 
$\int log{1 \over |u-z|} d\mu(u) + V_R(z)=0$ . From Theorem 5 we get
$$
\aligned
Re [ \int {{d\mu(u)} \over {u-z}} + V'(z) ] = Re [(R_{\mu}(z))^{1/2}]
\endaligned
\tag27
$$
where
$$
\aligned
R_{\mu}(z)= 
( V'(z))^2 - 2  \int_{supp\mu} {{V'(z)-V'(u)} \over {z-u}} d\mu(u)  \\ 
+ {1 \over z^2} (\int_{supp\mu} 2 (u+z) V'(u) ~ d\mu(u) ) . 
\endaligned
\tag28 
$$ 
This is a  function analytic in $K$, with possibly a pole
at zero. 
By integrating,  we have that $supp (\mu)$ is characterized by
$$
\aligned
Re \int^z (R_{\mu})^{1/2} dz =0.
\endaligned
\tag29
$$
The locus defined by (29) is a union of arcs with endpoints at zeros of
$R_{\mu}$.

Note that  
$$ 
\aligned  
R_{\mu} (z) \sim - [16 t^2 z^2 + \pi^2 (\rho^0 (z))^2] ,~~~~as~~~~z \to
\infty, \\
R_{\mu} (z) \sim {1 \over z^2} \int 2 u V'(u) d\mu (u), ~~~~as~~~~z \to 0.
\endaligned 
\tag30
$$

By conditions (1) for $\rho^0$,  $R_{\mu}$ is  blowing up  at the 
point $ \infty$ (at least for $t >0$; but the case 
$t=0$ is well understood: the equilibrium measure
consists of a single  analytic arc; see section 5). 
Hence it can only have finitely many zeros near infinity, 
otherwise they would have to accumulate near  
$\infty$ and then $R_{\mu}$ would be $0$ there.
On the other hand, $ R_{\mu} $ cannot have an accumulation point of zeros at
$z=0$,
because even if the pole at $0$ were removed (the coefficients of
${1 \over z^2}, ~{1 \over z}$ being zero),
$R_{\mu}$ would be holomorphically extended across $z=0$.
So, $R_{\mu}$ can only have a finite number of zeros in $\bar \Bbb K$.
It follows that
the support of the maximizing equilibrium measure consists of
only  finitely many arcs.

\bigskip

REMARK. Of course, conditions (1)  can be weakened. We could allow 
$\rho^0$ to have a pole at infinity of order other than two.
But our aim  here is not to prove the most general theorem possible, but 
instead illustrate a method that can be applied in the most general
settings under appropriate amendments.

\bigskip

REMARK. The assumption that $\rho^0$ 
is continuous and hence bounded at infinity  
is only needed to prove the finiteness of the components
of the support of the equilibrium measure of
the maximizing continuum.
If it is dropped then we may have an infinite number of components for 
isolated values of $x, t$. 
This will result in infinite genus representations of the semiclassical
asymptotics.
Of course infinite genus solutions of the focusing NLS equation
are known and well understood. So the analysis of [3] is expected to also
apply in that case, although it will be more tedious.

\bigskip

For a justification of the  "finite gap ansatz",
concerning the semiclassical limit of focusing NLS,
it only remains to verify the  "S-property".

\bigskip

7. THE S-PROPERTY

\bigskip

THEOREM 7. (The S-property)

Let $C$ be the contour maximizing the 
equilibrium energy, for the field given by 
(3a) with conditions (1).
Let $\mu$ be the extension of its equilibrium measure to the
full complex plane via (21a). 
Assume for simplicity that assumption (A)   holds.
Let $X(z) = \int_{supp\mu} log ({{ 1} \over {u-z}}) d\mu(u)$,
$X_R(z) = ReX(z) = \int_{supp\mu} log ({1 \over {|u-z|}}) d\mu(u), W_{\mu} =
X'$. 
Then, at any interior point of $supp\mu$ other than zero,
$$
\aligned
{d \over {d n_+}} (V_R + X_R) =
{d \over {d n_-}} (V_R + X_R),
\endaligned
\tag8a
$$ 
where the two derivatives above denote the normal derivatives,
on the $+$ and $-$ sides respectively.

PROOF: From Theorem 5, we have
$$
| Re (W_{\mu}(z) + V'(z) )| = | Re (R_{\mu})|^{1/2} .
$$

Using the definition for $X$, the above relation becomes
$$
|{d \over {dz}} Re(X+V) | =  | Re (R_{\mu})|^{1/2} .
$$

Now, $Re (X+V) =0 $ on the support of the
equilibrium measure. So, in particular $Re (X+V)$ is constant along the
equilibrium measure. Hence $|{d \over {dz}} Re (X+V)|$ must be equal to
the modulus of  $each$ normal derivative across the equilibrium measure.
So, 
$$
\aligned
|{d \over {d n_{\pm}}} (V_R + Re X)| = 
| {d \over {dz}} Re (X+V)| =  | Re (R_{\mu})|^{1/2} .
\endaligned
$$
Hence,
$$
\aligned
|{ d \over {d n_+}} (V_R + Re X) | =
|{d \over {d n_-}} (V_R + Re X)|.
\endaligned
$$ 
But it is easy  to see that both LHS and RHS quantities inside the
modulus sign are negative. This is because $V_R + ReX =0$ on
$supp\mu$ and negative on each side of $supp\mu$.
Hence result.

\bigskip

REMARK.
Once Theorem 7 is proved it follows by the Cauchy-Riemann equations that
$(V_I + Im X)_+ + (V_I + Im X)_-$ is constant
on each connected component of $supp\mu$, which means
that $Im \tilde \phi$ is constant on connected components of the contour,
where $\tilde \phi$ is as defined in formula (4.13) of [3].
This proves the existence of the appropriate "g-functions" in [3].

\bigskip

We recapitulate our results in the folllowing theorem, set
in the upper complex half-plane. Note that (8a) is the
"doubled up" version of  (8).

\bigskip

THEOREM 8. Let $\phi $ be given by (3a), 
where $\rho^0$ satisfies conditions (1).
Under assumption (A),
there is a piecewise smooth contour $C \in \Bbb F$, containing points $0_+, 0_-$ and 
otherwise lying in the cut upper half-plane $\Bbb K$,
with equilibrium measure $\lambda^C $,
such that $supp (\lambda^C )$ 
consists of a union of finitely many analytic arcs and
$$
\aligned
E_{\phi} (\lambda^C )  = max_{C' \in \Bbb F} E_{\phi} (\lambda^{C'} )
= max_{C' \in \Bbb F} [ inf_{\mu \in M(F)} E_{\phi} (\mu) ].
\endaligned
$$
On each interior point of $supp (\lambda^C )$ we have
$$
\aligned
{d \over {d n_+}} (\phi + V^{\lambda^C}) =
{d \over {d n_-}}  (\phi + V^{\lambda^C}),
\endaligned
\tag8
$$
where $V^{\lambda^C}$ is the  Green's potential of
the equilibrium measure $\lambda^C$ (see (5))
and the two derivatives above are the normal derivatives.

A curve satisfying (8) such that the support of its equilibrium measure
consists of a union of finitely many analytic arcs is called an S-curve.

PROOF: The fact that the  maximizing 
continuum $C$ is actually a contour is proved as follows.
If this were not the case, then we could choose
a subset of $C$, say $F$, which $is$ a contour, starting at $0_+$ and ending 
at $0_-$,
and going around the point $iA$. Clearly, by definition, the equilibrium energy
of $C$ is less than
the equilibrium energy of $F$, i.e. $E_{\phi} (\lambda^C)  \leq E_{\phi}
(\lambda^F) .$
On the other hand, since $C$ maximizes the equilibrium energy, we have
 $E_{\phi} (\lambda^F)  \leq E_{\phi} (\lambda^C) .$ So
 $E_{\phi} (\lambda^F)  = E_{\phi} (\lambda^C) .$

\newpage

8. CONCLUSION. 

\bigskip

In view of the interpretation of the variational problem in terms of the
semiclassical NLS problem, we have  the following result.

Consider the semiclassical limit ($\hbar  \to 0$)
of the solution of (9)-(10) 
with bell-shaped initial data. Replace the initial data by the so-called
soliton ensembles data (as introduced in [3])  defined by replacing
the scattering data for $\psi (x,0)= \psi_0(x)$ by their
WKB-approximation, so that
the spectral density of eigenvalues is
$$
\aligned
d\mu_0^{WKB}(\eta):=\rho^0(\eta)\chi_{[0,iA]}(\eta)d\eta +
\rho^0(\eta^*)^*\chi_{[-iA,0]}(\eta)d\eta, \\with~~~
\rho^0(\eta):= \frac{\eta}{\pi}\int_{x_-(\eta)}^{x_+(\eta)}\frac{dx}{
\sqrt{A(x)^2 +\eta^2}}
=
\frac{1}{\pi}\frac{d}{d\eta}\int_{x_-(\eta)}^{x_+(\eta)}\sqrt{A(x)^2+\eta^2}\,
dx,
\endaligned
$$
for $\eta \in (0,iA)$, where $x_- (\eta) < x_+ (\eta) $ are the two real turning
points, i.e.  $(A(x_{\pm}))^2 + \eta^2 =0$,
the square root is positive  
and the imaginary segments
$(-iA,0)$ and $(0,iA)$ are both considered to be oriented from bottom
to top to define the differential $d\eta$.

Assume that $\rho^0 $  
satisfies conditions (1).
Then, under assumption (A), asymptotically as $\hbar \to 0$, the
solution  $\psi (x,t)$ admits a "finite genus description".
(For a more precise explanation, see Appendix A2.)

The proof of this  is the main result  of [3], $assuming$ that the
variational problem of section 1 has an S-curve as a solution. But
this is now guaranteed by Theorem 8.

\bigskip

REMARK. For conditions weaker than the above,
the particular spectral density $\rho^0$ arising in
the semiclassical NLS problem can 
conceivably admit branch singularities
in the upper complex plane and condition (1) will not be satisfied. 
We claim  that even in such a case the finite 
gap genus  can be  justified, at least generically. The proof of this fact will
require setting the variational problem on a 
Riemann surface with moduli at the  branch singularities of  $\rho^0$.

\newpage

REMARK. Consider the semiclassical problem (9)-(10)
in the case of initial data $\psi_0(x) = Asechx$, where $A >0$. Then the
WKB density is given by  $\rho^0=i$ (see (3.1) and (3.2) of [3]; note that
condition (1)   is satified). 
So the finite genus ansatz holds for any $x,t$,
as long as  the assumption (A) of section 5 holds. But then assumption (A)
can be eventually dropped; see Appendix A3.

\bigskip

REMARK. The behavior of a solution of (9) in general depends not only on the
eigenvalues of the Lax operator, but also on the associated norming constants and
the reflection coefficient.
In the special case of the soliton ensembles data the norming constants alternate
between $-1$ and $1$ while the reflection coefficient is by definition zero.
More generally, for real analytic data  decaying at infinity the
reflection coefficient is exponentially small everywhere except at zero and
can be neglected (although the rigorous proof of this is not trivial).
 
\bigskip

ACKNOWLEDGEMENTS.

The first author acknowledges the kind support of the General Secretariat of
Research and Technology, Greece, in particular grant 97EL16. He is also grateful
to the
Department of Mathematics of the
University of South Florida for its hospitality during a visit on May 2001 and
to the Max Planck Society for support since 2002.
Both authors  acknowledge the invaluable contribution  of our collaborators
Ken McLaughlin and Peter Miller through  stimulating
discussions, important comments, corrections and constructive criticism.

\newpage

APPENDIX A1. COMPACTNESS OF  THE SET OF CONTINUA

\bigskip

In this section we prove that 
the sets $I(\bar \Bbb K)$ and hence
$\Bbb F$ defined in section 1 are compact and complete.

As stated in section 1, the space we must work with is the upper half-plane:
$ \Bbb H = \{ z: Im z >0 \} $.
The closure of this space  is $\bar \Bbb H =  \{ z: Im z \geq 0 \} \cup
\{\infty\}$.
Also
$ \Bbb K = \{ z: Im z >0 \} \setminus \{ z: Rez =0, 0< Im z \leq A \}$.
In the closure of this space, $\bar \Bbb K $, we consider the points
$ix_+$ and $ix_-$, where $0 \leq x < A$ as distinct.

Even though we eventually wish to consider only smooth contours, we are forced
to a
priori work with general closed sets. The reason is that the set of contours is
not compact in any
reasonable way, so it seems impossible to prove any existence theorem for a
variational problem
defined only on contours. Instead,
we define $\Bbb F$ to be the set of all "continua"
$F$ in $\bar \Bbb K$
(i.e. connected compact sets, containing the points $0_+, 0_-$).

Furthermore, we need to introduce an appropriate topology on $\Bbb F$, that will
make it
a compact set. In this we follow the discussion of  Dieudonn\'e 
([6], chapter III.16).

We think of the closed upper half-plane $\bar \Bbb H$ as a compact space in the 
Riemann sphere. We thus choose to equip
$\bar \Bbb H$ with the  "chordal" distance, denoted by $d_0 (z, \zeta)$, 
that is  the distance between
the images of $z$ and $\zeta$ under the stereographic projection.
This induces naturally a distance in $\bar \Bbb K $ (so, for example,
$d_0 (0_+, 0_-) \neq 0).$
We also denote by $d_0$ the induced  distance between compact sets
$E, F$ in $\bar \Bbb K$:
$d_0(E,F) = max_{z \in E} min_{\zeta \in F} d_0 (z,\zeta)$. 
Then, we define the so-called Hausdorff metric on the set $ I ( \bar \Bbb K   )$
of closed nonempty subsets of  $\bar \Bbb K$
as follows. 
$$
\aligned
d_{\Bbb K} (A,B) = sup  ( d_0 (A,B), d_0 (B,A) ).
\endaligned
\tag A.1
$$

LEMMA A.1. The Hausdorff metric defined by (A.1) is indeed a metric.
The set $ I ( \bar \Bbb K   )$
is compact and complete.

PROOF: It is clear that $d_{\Bbb K} (A,B)$ is non-negative
and symmetric by definition.
Also if $d_{\Bbb K} (A,B) =0 $, then $ d_0 (A,B)=0$, hence
$ max_{z \in A} min_{\zeta \in B} d_0 (z,\zeta) =0$ and thus for all
$z \in A$, we have $ min_{\zeta \in B} d_0 (z,\zeta) =0$. In other words,
$z \in B$. By symmetry, $A=B$.

The triangle inequality follows from the triangle inequality for $ d_0 $. 
Indeed,  suppose $A, B, C \in   I ( \bar \Bbb K   )$. Then
$d_{\Bbb K} (A,B) =  sup  (d_0 (A,B), d_0 (B,A) )  =
d_0 (A,B) ,$ without loss of generality. Now,
$$ 
\aligned
d_0 (A,B) = max_{z \in A} min_{\zeta \in B} d_0 (z,\zeta) 
\\
\leq   max_{z \in A} min_{\zeta \in B}  min_{\zeta_0 \in C}
(d_0 (z,\zeta_0) +  d_0 (\zeta_0, \zeta)) ,
\endaligned
$$
by the
triangle inequality for $ d_0 $. Let $z=z_0  \in A$ be the value of $z$ 
that maximizes
$ min_{\zeta \in B}   min_{\zeta_0 \in C}
(d_0 (z,\zeta_0) + d_0 (\zeta_0, \zeta)).$
This is then
$$ 
\aligned
min_{\zeta \in B}  min_{\zeta_0 \in C}
(d_0 (z_0,\zeta_0) + d_0 (\zeta_0, \zeta))  \\
\leq  min_{\zeta_0 \in C} d_0 (z_0,\zeta_0)  +
 min_{\zeta \in B} min_{\zeta_0 \in C} d_0 (\zeta_0, \zeta)  \\
\leq max_{z \in A}  min_{\zeta_0 \in C}
d_0 (z,\zeta_0) + max_{\zeta \in B}  min_{\zeta_0 \in C} d_0 (\zeta_0, \zeta) \\
\leq d_0 (A,C) + d_0 (B,C) \leq d_{\Bbb K} (A,C) + d_{\Bbb K} (B,C). 
\endaligned
$$
The result follows from symmetry.

We will next show that $ I ( \bar \Bbb K   )$ is complete
and precompact. Since a precompact, complete metric space is compact
([6], proposition (3.16.1)) the proof of Lemma A.1 follows.

LEMMA A.2  If the metric space $ \Bbb E $ 
equipped with a distance $d_0$ is complete, then so is  $ I (  \Bbb E   )$,
the set of closed nonempty subsets of $ \Bbb E $ , equipped with the
Hausdorff distance
$$
\aligned
d_{\Bbb E} (A,B) = sup  ( d_0 (A,B), d_0 (B,A) ),
\endaligned
$$
for any closed nonempty subsets 
$A, B$, where $d_0 (A,B) = max_{a \in A} min_{b \in B} d_0 (a,b)$.

Furthermore,
if  $ \Bbb E $ is precompact, then so is  $ I (  \Bbb E   )$.

PROOF: Suppose $ \Bbb E $ is complete.
Let $X_n$ be a Cauchy sequence in  $ I ( \Bbb E   )$. We will
show that  $X_n$ converges to $X= \cap_{n \geq 0} \bar \cup_{p \geq 0} X_{n+p}.$
(Overbar denotes closure.) 

Indeed, given any $\epsilon >0$, 
$$
\aligned
d_0 (X_n, X) = 
max_{x \in X_n} min_{y \in X} d_0 (x,y)\\
\leq max_{x \in X_n} max_{y \in \bar \cup_{p \geq 0}  X_{n+p}}  d_0 (x,y)
<\epsilon ,
\endaligned
$$
for large $n$, by the completeness of $ \Bbb E $.
Similarly,
$$
\aligned
d_0 (X, X_n) = 
max_{x \in X} min_{y \in X_n} d_0 (x,y)\\
\leq max_{x \in  \bar \cup_{p \geq 0} X_{n+p}} min_{y \in  X_{n}}  d_0 (x,y)
<\epsilon.
\endaligned
$$

Next,  suppose $ \Bbb E $ is precompact. Then, by definition, 
given any $\epsilon >0$, there is a finite set, say 
$S=\{s_1, s_2, ..... , s_n \}$,
where $n$ is a finite integer, such that any point $x$ of
$ \Bbb E $ is at a distance $d_0$
less than $\epsilon$ to the set $S$.
Now,  consider the set  of  subsets of $S$,
 which is of course finite. Clearly every closed set is at a distance
less than $\epsilon$ to a member of that set: 
$$
\aligned
d_0  (A, S) = max_{a \in A} min_{s \in S}  d_0 (a,s) <\epsilon ,\\
d_0  (S, A)= max_{s \in S} min_{a \in A} d_0(a,s) <\epsilon,
\endaligned
$$
for any closed nonempty set $A$. 
Hence $d_{\Bbb E} (A, S)<\epsilon.$

So, any closed nonempty set $A$ is at a distance less than 
$\epsilon$ to the finite power set of $S$.
So $ I (  \Bbb E   )$ is precompact.

\bigskip

APPENDIX A2. THE DESCRIPTION OF THE SEMICLASSICAL LIMIT OF THE FOCUSING NLS
EQUATION
UNDER THE FINITE GENUS ANSATZ

\bigskip

We present one of the main results of [3] on 
the semiclassical asymptotics for problem (9)-(10),
in view of the fact that 
the finite genus ansatz holds.
In particular, we fix $x,t$ and use the result that the support of the maximizing
measure of
Theorems 4 and 8 consists of a finite union of analytic arcs.

First, we define the so-called g-function. 
Let $C$ be the  maximizing contour of Theorem 4.
A priori we seek a function satisfying
$$
\aligned
g(\lambda) ~is ~independent ~of ~ \hbar.\\
g(\lambda) ~is ~~analytic ~~~for
\lambda\in \Bbb C \setminus (C\cup C^*).\\
g(\lambda)\rightarrow 0 ~~as~~
\lambda\rightarrow\infty.\\
g(\lambda) ~~~assumes~~
continuous ~~boundary~~
values ~~from ~~both ~~~sides ~~~~of ~~~C\cup C^*,\\
denoted~~by~~~g_+(g_-)~~on~~the~~left~~(right)~~~of~~~C \cup C^*.\\
g(\lambda^*)+g(\lambda)^* = 0 ~~~~~for~~~
all ~\lambda\in \Bbb C\setminus (C\cup C^*).
\endaligned
$$

The assumptions above are satisfied if we write $g$ in terms of the 
maximizing equilibrium measure of Theorem 8,
$d\mu = d\lambda^C = \rho (\eta ) d\eta,$ doubled up according to (21a).  Indeed,
$$
\aligned 
g(\lambda) =  \int_{C \cup C^*} log(\lambda -\eta) \rho(\eta) d\eta,
\endaligned
$$
for an appropriate definition of the logarithm branch (see [3]).

For $\lambda\in C$, define
the functions
$$
\aligned 
\theta(\lambda):=i(g_+(\lambda)-
g_-(\lambda)),\\ 
\Phi(\lambda)  
:=
\int_{0}^{iA} log(\lambda-\eta)\rho^0(\eta)\,d\eta  +
\int_{-iA}^0 log(\lambda-\eta)\rho^0(\eta^*)^*\,d\eta \\
+
2i\lambda x + 2i\lambda^2 t + i\pi \int_\lambda^{iA}\rho^0(\eta)\,d\eta -
g_+(\lambda) - g_-(\lambda),
\endaligned
$$
where $\rho^0(\eta)$ is the holomorphic function
(WKB density of eigenvalues) introduced in section 1 
(see conditions (1)).

The finite genus ansatz implies that for each $x, t$ there is a finite positive
integer
$G$ such that the contour $C$ can be divided into "bands" 
[the support of $\rho(\eta) d\eta$] and "gaps" (where $\rho =0$).
We denote these bands by $I_j$. More precisely,
we define the   analytic arcs  $I_j, I^*_j,  j=1, ... , G/2$ as follows (they
come in conjugate pairs).
Let the points  $\lambda_j,~~~j=0, ... , G$, in the open upper half-plane
be the branch points of the function $g$.
All such points lie
on the contour $C$ and we order them as
$\lambda_0, \lambda_1, ... , \lambda_G$,
according to the direction given to $C$.
The points $\lambda_0^*, \lambda_1^*, ... , \lambda_G^*$ are their
complex conjugates.
Then let
$I_0 = [0, \lambda_0]$ be the subarc of $C$ joining points $0$ and $\lambda_0$.
Similarly,
$ I_j = [\lambda _{2j-1}, \lambda _{2j}],~~~ j=1, ... , G/2$.
The connected components of the set $\Bbb C \setminus \cup_j (I_j \cup I^*_j)$
are the so-called "gaps", for example the gap $\Gamma_1$ joins $\lambda_0$
to $\lambda_1$, etc.

It actually follows from the properties of
$g, \rho$ that 
the function $\theta (\lambda)$ defined on $C$ is constant on each of the
gaps $\Gamma_j$, taking a value which we will denote by $\theta_j$, while the 
function $\Phi$ is constant on each of the bands,
taking the value denoted by $\alpha_j$ on the band $I_j$.

The finite genus ansatz for the given fixed $x, t$ implies 
that the asymptotics of the solution of (9)-(10)
as $\hbar  \to 0$ can be given by the next theorem.

\bigskip

FINITE GAP ANSATZ THEOREM A.1.
Let $x_0, t_0$ be given.
The solution $\psi(x,t)$ of (9)-(10) is asymptotically
described  (locally) as a slowly
modulated $G+1$ phase wavetrain.  Setting $x=x_0+ \hbar  \hat{x}$ 
and $t=t_0+\hbar  \hat{t}$,
so that $x_0, t_0$ are "slow" variables
while $\hat{x}, \hat{t}$ are "fast" variables,
there exist
parameters

$a,  U = (U_0, U_1, .... , U_G)^T,$
$k =(k_0, k_1, ......, k_G)^T,$
$w =(w_0, w_1, ....., w_G)^T, $
$Y =(Y_0, Y_1, ........., Y_G)^T,$
$Z =( Z_0, Z_1, ...... , Z_G)^T $
depending on the slow variables
$x_0$ and $t_0$  (but not  $\hat{x}, \hat{t}$)
such that
$$
\aligned
\psi(x,t) = \psi (x_0 + \hbar  \hat{x},  t_0 + \hbar  \hat{t}) \sim
a(x_0, t_0) e^{iU_0(x_0, t_0)/\hbar}
e^{i(k_0(x_0, t_0) \hat{x}-w_0(x_0, t_0) \hat{t})} \\
\cdot  \frac{\Theta(  Y(x_0, t_0)+
i  U(x_0, t_0)/\hbar +
i(  k(x_0, t_0) \hat{x}-  w(x_0, t_0)\hat{t}))}
{
\Theta(  Z(x_0, t_0)+
i  U(x_0, t_0)/\hbar +
i(  k(x_0, t_0)\hat{x}-  w(x_0, t_0)\hat{t}))}.
\endaligned
\tag A.2
$$

All parameters can be defined in terms of an underlying
Riemann surface $X$.
The moduli of $X$ are given by $\lambda_j, ~ j=0, .... , G$ and their complex
conjugates  $\lambda_j^*, ~ j=0, .... , G$.
The genus of $X$ is $G$. The moduli of $X$ vary slowly with $x,t$, i.e.
they depend on  $x_0, t_0$ but not
$\hat{x}, \hat{t}$. For the exact formulae 
for the parameters
as well as the definition  of the theta functions we present the following
construction.
 
The Riemann surface $X$ is constructed by cutting two copies of the complex
sphere along the slits
$I_0 \cup I_0^*, I_j, I_j^*, j=1. ... ,G$,
and pasting the "top" copy to the "bottom" copy along these
very slits.

We define the homology cycles
$a_j, b_j, ~~~ j=1, ... , G$
as follows.
Cycle $a_1$ goes around
the slit $I_0 \cup I_0^*$ joining $\lambda_0$ to
$\lambda_0^*$,
remaining on the top sheet, oriented counterclockwise,
$a_2$ goes through  the slits $I_{-1}$ and $I_1$
starting from the top sheet, also
oriented counterclockwise,
$a_3$ goes around the slits $I_{-1}, I_0 \cup I_0^*, I_1$
remaining on the top sheet, oriented counterclockwise, etc.
Cycle $b_1$ goes through $I_0$  and $I_1$ oriented counterclockwise,
cycle  $b_2$ goes through $I_{-1}$  and $I_1$, also
oriented counterclockwise,
cycle  $b_3$ goes through $I_{-1}$  and $I_2$, and
around the slits $I_{-1}, I_0 \cup I_0^*, I_1$,
oriented counterclockwise, etc.

On $X$
there is a complex $G$-dimensional linear space of holomorphic
differentials, with basis elements $\nu_k(P)$ for $k=1,\dots,G$ that
can be written in the form
$$
\aligned
\nu_k(P)=\frac{\displaystyle\sum_{j=0}^{G-1}c_{kj}\lambda(P)^j}
{R_X(P)}\,d\lambda(P)\,,
\endaligned
$$
where $R_X(P)$ is a ``lifting'' of the function $R(\lambda)$
from the cut plane to $X$: if $P$ is on the first sheet of $X$ then
$R_X(P)=R(\lambda(P))$ and if $P$ is on the second sheet of
$X$ then $R_X(P)=-R(\lambda(P))$.  The coefficients $c_{kj}$
are uniquely determined by the constraint that the differentials
satisfy the normalization conditions:
$$
\aligned
\oint_{a_j}\nu_k(P) = 2\pi i\delta_{jk}.
\endaligned
$$
From the normalized differentials, one defines a $G\times G$ matrix
$H$ (the period matrix) by the formula:
$$
\aligned
H_{jk}=\oint_{b_j}\nu_k(P).
\endaligned
$$
It is a consequence of the standard theory of Riemann surfaces that
$ H$ is a symmetric matrix whose real part is  negative
definite.
 
In particular, we can
define the theta function
$$
\aligned
\Theta( w):=\sum_{  n\in {\Bbb Z}^G}\exp({1 \over 2} n^T  H
n+ n^T  w),
\endaligned
$$
where $ H$ is the period matrix associated to $X$.
Since the real part of $ H$ is  negative
definite, the series converges.
 
We arbitrarily fix a base point $P_0$ on $X$.  The
Abel map
$ A:X\to Jac(X)$ is then defined componentwise as follows:
$$
\aligned
A_k(P;P_0):=\int_{P_0}^P\nu_k(P'),~~~~ k=1,\dots,G,
\endaligned
$$
where $P'$ is an integration variable.
 
A particularly important element of the Jacobian is the
Riemann constant vector
$K$ which is defined, modulo the lattice
$\Lambda$, componentwise by
$$
\aligned
K_k:=\pi i + \frac{H_{kk}}{2}-\frac{1}{2\pi i}\sum_{j=1\atop j\neq
k}^G\oint_{a_j} \left(\nu_j(P)\int_{P_0}^P\nu_k(P')\right),
\endaligned
$$
where the index $k$ varies between $1$ and $G$.

Next, we will need to define a certain meromorphic differential
on $X$.  Let $\Omega(P)$ be holomorphic away from the points $\infty_1$
and $\infty_2$, where it has the behavior
$$
\aligned
\Omega(P) = dp(\lambda(P)) + \displaystyle
\left(\frac{d\lambda(P)}{\lambda(P)^2}\right), ~~~P\to
\infty_1,\\
\Omega(P) =-dp(\lambda(P))+\displaystyle
O\left( \frac{d\lambda(P)}{\lambda(P)^2}\right), ~~~~
P\to \infty_2,
\endaligned
$$
and made unique by the normalization conditions
$$
\aligned
\oint_{a_j}\Omega(P) = 0, j=1,\dots,G.
\endaligned
$$
Here $p$ is a  polynomial, defined as follows.

First, let us  introduce the function $R(\lambda)$
defined by
$$
\aligned
R(\lambda)^2 = \prod_{k=0}^{G}(\lambda-\lambda_k)(\lambda-\lambda_k^*),
\endaligned
$$
choosing the particular branch that is cut along the bands $I_k^+$ and
$I_k^-$ and satisfies
$$
\aligned
\lim_{\lambda\rightarrow\infty}\frac{R(\lambda)}{\lambda^{G+1}}= -1.
\endaligned
$$
This defines a real function,  i.e. one that satisfies
$R(\lambda^*)=R(\lambda)^*$.  At the bands, we have
$R_+(\lambda)=-R_-(\lambda)$, while $R(\lambda)$ is analytic in the
gaps.  Next, let us introduce the function $k(\lambda)$ defined by
$$
\aligned
k(\lambda)=\frac{1}{2\pi i}
\sum_{n=1}^{G/2}\theta_n\int_{\Gamma_n^+\cup\Gamma_n^-}
\frac{d\eta}{(\lambda-\eta)R(\eta)} +
\frac{1}{2\pi i}\sum_{n=0}^{G/2}
\int_{I_n^+\cup I_n^-}\frac{\alpha_n 
~~d\eta}{(\lambda-\eta)R_+(\eta)}\,.
\endaligned
$$
Next  let
$$
\aligned
H(\lambda)=k(\lambda)R(\lambda).
\endaligned
$$
The function $k$ satisfies
the jump relations
$$
\aligned
k_+(\lambda)-k_-(\lambda)=
\displaystyle -\frac{\theta_n}{R(\lambda)},~~~~\lambda\in
\Gamma_n^+\cup\Gamma_n^-\\
k_+(\lambda)-k_-(\lambda)=
\displaystyle -\frac{\alpha_n}{R_+(\lambda)},~~~~\lambda\in
I_n^+\cup I_n^-,
\endaligned
$$
and is otherwise analytic. It 
blows up like $(\lambda-\lambda_n)^{-1/2}$ near each
endpoint, has continuous boundary values in between the endpoints, and
vanishes like $1/\lambda$ for large $\lambda$.  It is the only such
solution of the jump relations.
The factor of $R(\lambda)$ renormalizes
the singularities at the endpoints, so that, as desired, the boundary
values of $H(\lambda)$ are bounded continuous functions.  Near
infinity, there is the asymptotic expansion:
$$
\aligned
H(\lambda)=H_G\lambda^G + H_{G-1}\lambda^{G-1} +
\dots + H_1\lambda + H_0 + O (\lambda^{-1})\\
=p(\lambda) + O (\lambda^{-1}),
\endaligned
\tag A.3
$$
where all coefficients $H_j$ of the polynomial $p(\lambda)$
can be found explicitly by expanding $R(\lambda)$ and the Cauchy
integral $k(\lambda)$ for large $\lambda$.  It is easy
to see from the reality of $\theta_j$ and $\alpha_j$ that
$p(\lambda)$ is a polynomial with real coefficients.

Thus the polynomial $p(\lambda)$ is defined and hence the
meromorphic differential  $\Omega(P)$ is defined.

\bigskip

Let the vector $U\in {\Bbb C}^G$ be defined
componentwise by
$$
\aligned
U_j:=\oint_{b_j}\Omega(P).
\endaligned
$$
Note that $\Omega(P)$ has no residues.

Let the  vectors $ V_1,  V_2$ be defined componentwise by
$$
\aligned
V_{1,k}=
(A_k(\lambda_{1+}^*)+A_k(\lambda_{2+})+A_k
(\lambda_{3+}^*)+\dots+A_k(\lambda_{G+}))
+A_k(\infty) +\pi i + \frac{H_{kk}}{2},\\
V_{2,k}=(A_k(\lambda_{1+}^*)+A_k(\lambda_{2+})+A_k
(\lambda_{3+}^*)+\dots+A_k(\lambda_{G+}))
-A_k(\infty) +\pi i + \frac{H_{kk}}{2},
\endaligned
$$
where $k =1, ..., G$, and  the $+$ index means that the integral for $A$
is to be taken on the first sheet of $X$, with base point $\lambda^0_+.$

Finally, let
$$
\aligned
a=\frac{\Theta( Z)}
{\Theta(  Y)}
\sum_{k=0}^G(-1)^k\Im(\lambda_k)  ,\\
k_n=\partial_x U_n,~~~~~
w_n=-\partial_t U_n,~~~~~
n=0,\dots,G,
\endaligned
$$
where
$$
\aligned 
Y=- A(\infty)-  V_1,~~~~ 
Z=  A(\infty)- V_1,
\endaligned
$$
and  $ U_0 =-(\theta_1 + \alpha_0)$
where $\theta_1$ is the (constant in $\lambda$) value of the function $\theta$
in the gap $\Gamma_1$
and $\alpha_0$ is the (constant) value of the function $\phi$ in the band $I_0$.

Now, the parameters appearing in formula (A.2) are completely described.
 
We simply note  here that the $U_i$ and hence the $k_i$ and $w_i$ are real.  
We also note that the denominator in (A.2) never vanishes (for any
$x_0, t_0, \hat{x}, \hat{t}$).
 
\bigskip

REMARK. The most general version of Theorem A.1 is not fully proved in  this paper.
So far the main text of this paper and the analysis of [3] provide a proof under assumption (A).
Theorem A.2 is more general, because assumption (A) is dropped.
Appendix A4  shows how to remove the assumption of existence of an analytic extension of the
limiting density of eigenvalues. 
But there is a remaining issue: the validity of the solitons ensemble approximation. This final question 
can be answered via the so-called exact WKB theory; a related publication (with Setsuro Fujiie) is forthcoming.
\bigskip

REMARK. Theorem A.1 presents pointwise asymptotics in $x, t$. In [3], these are
extended
to uniform asymptotics in certain compact sets covering the $x,t$-plane.
Error estimates are also given in [3].

\bigskip

REMARK. As mentioned above, we do not know if the support of the equilibrium measure
of the maximizing continuum is unique. But the asymptotic formula (A.2) depends
only on the endpoints $\lambda_j$ of the analytic subarcs of the support.
Since the asymptotic expression (A.2) must be unique, it is easy to see that 
the endpoints also must be unique. Different Riemann surfaces give different 
formulae (except of course in degenerate cases: a degenerate genus 2 surface 
can be a pinched genus 0 surface and so on).

\bigskip

APPENDIX A3. DROPPING ASSUMPTION (A) OF SECTION 5.

\bigskip

This appendix is presented as appeared in a corrected form in
the Journal of Mathematical Physics, v.50, n.9, 2009, signed by one of us (S.K.).  
\bigskip

In section 5, we have assumed that the solution of the problem of the
maximization
of the equilibrium energy  is a continuum, say $F$, which does not intersect the
linear
segment $[0, iA]$ except of course  at $0_+, 0_-$.
We also prove that $F$ does not touch the real line, except of course at
$0$ and possibly $\infty$. This enables us to take variations in section 6 of [9],
keeping fixed a finite number of points, and thus arrive at the identity of
Theorem 5,
from which we derive the regularity of $F$ and the fact that $F$ is, after all,
an S-curve.

In general, it is conceivable that $F$ intersects the linear
segment $[0, iA]$ at points other than  $0_+, 0_-$.
If the set of such points is finite, there is no problem, since we can always
consider
variations keeping fixed a finite number of points, and arrive at the same
result (see the remark after the proof of Theorem 5).

If, on the other hand, this is not the case, we have a different kind of
problem,
because the function $V$ introduced in section 6 (the complexification of
the field) is not analytic across the segment $[-iA, iA]$.

What is true, however, is that $V$ is analytic in a Riemann surface consisting
of infinitely
many sheets, cut along the line segment $[-iA, iA]$. So, the appropriate,
underlying space for the
(doubled up) variational problem should now be a
non-compact Riemann surface, say $\Bbb L$.

Compactness is crucial in the proof of a maximizing continuum. But we can
compactify the
Riemann surface $\Bbb L$ by compactifying the
complex plane.
Let  the map $\Bbb C \to \Bbb L$ be defined by
$$
\aligned
y = log(z-iA) - log(z+iA).
\endaligned
$$
The point $z =iA$ corresponds to infinitely many y-points, i.e.
$y = -\infty + i \theta,~~\theta \in \Bbb R$, which will be identified.
Similarly,
the point $z= -iA$  corresponds to infinitely many points
$y = +\infty + i \theta,~~\theta \in \Bbb R$, which will also be identified.
The point $0 \in \Bbb C$  corresponds to the points $ k \pi i$, $k$ odd.
 
By compactifying the plane we then compactify the Riemann surface $\Bbb L$.
The distance between two points in the Riemann surface $\Bbb L$
is defined to be the corresponding stereographic distance between the
images of these points in the compactified $\Bbb C$.

With these changes, the proof of the existence of the maximizing continuum in
sections 1, 3, 4 
goes through virtually unaltered. In section 6, we would have to
consider the complex field $V$ as a function defined in the
Riemann surface $\Bbb L$ and all proofs go through.
The corresponding
result of section 7 
will give us an S-curve $C$ in the Riemann surface $\Bbb L$.
We then have the following facts.

Consider the  image $\Bbb D$ of the closed upper half-plane
under
$$
\aligned
y = log(z-iA) - log(z+iA).
\endaligned
$$
Consider continua in $\Bbb D$ containing the points $y=\pi i$ and $y=-\pi i$.
Define the Green's potential  and Green's energy of a Borel measure by
(4), (5), (6)  and the equilibrium measure by (7).
Then there exists a continuum
$F$  maximizing the equilibrium energy, for the field given by
(3)  with conditions (1). $F$ does not touch $\partial \Bbb D$ except at a finite
number of points. By taking  variations as in section 6, one sees
that $F$ is an S-curve. In particular,
the support of the equilibrium measure on $F$ is a union of analytic
arcs and
at any interior point of $supp\mu$
$$
\aligned
{d \over {d n_+}} (\phi + V^{\lambda^F}) =
{d \over {d n_-}}  (\phi + V^{\lambda^F}),
\endaligned
$$
where the two derivatives above denote the normal derivatives.

\bigskip

We then have the following.

THEOREM A.2. Consider the semiclassical limit ($\hbar \to 0$)
of the solution of (9)-(10) (that is the initial value problem for the
focusing NLS with  parameter $\hbar$)
with bell-shaped initial data. Replace the initial data by the so-called
soliton ensembles data (as introduced in [3])  defined by replacing
the scattering data for $\psi (x,0)= \psi_0(x)$ by their
WKB-approximation.  Assume, for simplicity,
that the spectral density of eigenvalues
satisfies conditions (1).

Then, asymptotically as $\hbar \to 0$, the
solution  $\psi (x,t)$ admits a "finite genus description", in the sense of
Theorem A.1.

PROOF: (i) The proof of the existence of an S-curve $F$ in $\Bbb L$ follows as above.
It consists of a finite number of bands (the components of the
support of the equilibrium measure) and gaps.

(ii) We want to deform the original discrete Riemann-Hilbert problem to the set
$\hat F$ consisting of the  projection of  $F$
to the complex plane.
It is clear however that  $\hat F$ may not encircle the spike $[0,iA]$.
It is possible, on the other hand,  to append S-loops (considered in  $\Bbb L$)
and end up with a sum of S-loops, such that the amended $\hat F$  $does$
encircle the spike $[0,iA]$, meaning that  $[0, iA]$ is a subset
of the closure of the union   of the interiors of the loops
of which $ \hat F$ consists. A little thought shows that this is all we need.
(Indeed, within each of the loops we use the same pole-removing transformation as
in [3]. Eventually of course we have to  use different interpolations, according to the sheet
of each piece of $F$.)

To see that we can always append the needed S-loop,  suppose
there is an open interval, say $(i\alpha, i\alpha_1)$,
which lies in the exterior of $\hat F$, while
$i\alpha, i\alpha_1 \in \hat F$. Let us assume  for example that $\hat F$ crosses
$[0,iA]$ along bands at $i\alpha, i\alpha_1$ (these bands, say  $S, S_1,$ actually
belong to $F$ to be more precise) and also assume without loss of generality  that they
both lie in
the principal sheet.
Let $\beta^-, \beta^+$ be points (considered in $\Bbb C$)
lying on $S$ to the left and right
of $i\alpha$ respectively, and at a small distance
from $i\alpha$. Similarly,
let $\beta_1^-, \beta_1^+$ be points lying on $S_{1}$
to the left and right
of $i\alpha_{1}$ respectively, and at a small distance
from $i\alpha_{1}$.
We will show that there exists a "gap" region including the preimages of
$\beta^-, \beta_{1}^-$ lying in the $N$th sheet for $-N$ large enough,
and similarly there exists a "gap" region including the preimages of
$\beta^+, \beta_{1}^+$ lying in the $M$th sheet for $M$ large enough,
both being  regions for which the gap inequalities hold a priori,
irrespectively of the actual S-curve, depending only on the external field!

Indeed, note  that the
quantity  $Re (\tilde \phi^{\sigma} (z) )$ (which defines the variational
inequalities) is a priori bounded above  by
$-\phi(z)$. For this, see (8.8) in Chapter 8 of [3]; there is actually a sign error:
the right formula is
$$
\aligned
Re (\tilde \phi^{\sigma} (z) )= -\phi(z)
+ \int G(z,\eta) \rho^{\sigma} (\eta)d\eta.
\endaligned
$$
Next note (see for example (5.8) of [3] with $K$ varying along the natural numbers according
to the relevant sheet of $\Bbb L$)
that the difference of the values of the
function $Re (\tilde \phi^{\sigma} (z) )$ in consecutive sheets is
$\delta Re (\tilde \phi^{\sigma})  \sim \pm 2 Im \rho (z) \pi Rez$ near the spike
$[0, iA]$ (remember $Im \rho (z) >0$ there) and hence the difference of
the values at points on consecutive sheets whose image under the
projection to the complex plane is $i \eta + \epsilon$, where $\eta $ is real
and $\epsilon $ is a small (negative or positive) real,
is $\delta (Re \tilde \phi^{\sigma} )  \sim \pm 2 \pi Im \rho (z) \epsilon$.
This means that on the left (respectively right)
side of the imaginary semiaxis, the inequality $Re ( \tilde \phi^{\sigma} (z) )<0$
will be eventually (depending
on the sheet) be valid at any given small distance to it.

We now connect  the preimages of $\beta^-$
and  $\beta_{1}^-$ (under the projection of $\Bbb L $ to
$\Bbb C$) lying in the $N$th sheet to the preimages of
$\beta^-$ and $\beta_1^-$ lying in the principal sheet respectively. Similarly  we join
the preimages of $\beta^+$
and  $\beta_{1}^+$ lying in the $M$th sheet
to the preimages of
$\beta^-$ and $\beta_1^-$ lying in the principal sheet respectively.

Then, we  join the  the preimages of $\beta^-$
and  $\beta_{1}^-$ (under the projection of $\Bbb L $ to
$\Bbb C$) lying in the $N$th sheet and the preimages of $\beta^+$
and  $\beta_{1}^+$ lying in the $M$th sheet, along the
according gap regions.

It is easy to see that (together with the bands $S$ and $S_1$) we  end up with an S-loop
(in $\Bbb L$) whose projection is
covering the "lacuna" $(i\alpha, i\alpha_{1})$.

The original discrete Riemann-Hilbert problem  can be  trivially deformed to a
discrete Riemann-Hilbert on the resulting  (projection of the)
union of S-loops. All this is possible even in the
case where $\hat F$ self-intersects.

(iii) We  deform the discrete Riemann-Hilbert problem to the continuous one with the right
band/gap structure
(on $\hat F$; according to the  equilibrium measure on $F$), which is
then explicitly solvable via theta functions exactly as in [3].
Both the discrete-to-continuous approximation and the opening of the lenses needed
for this deformation are justified  as in [3]
and therefore the technical details will not be repeated here.
It is important to notice that our construction has ensured the analytic continuation of
the jump matrix along $\hat F$ (oriented according to $F$).
The g-function is defined by the same Thouless-type formula with respect to the
equilibrium measure (cf. section 2(iii)).
It satisfies the same conditions as in [3]
(measure reality and variational inequality) on bands  and  gaps.
The equilibrium measure lives in $\Bbb L$ but the  Riemann-Hilbert problem lives in $\Bbb C$.

\newpage

APPENDIX A4.  DROPPING THE ASSUMPTION OF AN ANALYTIC EXTENSION OF THE SPECTRAL DENSITY  $\rho^0 $.  

\bigskip

This section has previously appeared as a Max Planck Institute preprint in 2002, signed by one of us (S.K.).

\bigskip

THEOREM A.3. The finite gap ansatz Theorem A.1 is valid for the solution of the problem
(9)-(10), if we substitute the initial data by their soliton ensembles approximation.

No assumption of an analytic extension for  $\rho^0 $ is necessary.
\bigskip

SKETCH OF PROOF:

It is essential for the proofs in [3] that the "density
of eigenvalues" $\rho^0(\eta)$ (see (3.2) of [3]),
derived by WKB theory and a priori
defined in the straight line interval connecting
$0$ to $iA$, be analytically extensible
to the closed upper half-plane  $\Bbb H$.
The main issue is whether the function
$$
\aligned
R^0(\eta) = \int^{x_+(\eta)}_{x_-(\eta)} (A(x)^2 +\eta^2)^{1/2} dx,
\endaligned
$$
where the turning points are defined by
$$
\aligned
A( x_{\pm} (\eta)) = -i\eta, ~~~0< -i\eta < A, \\
-A < x_- (\eta) < 0 < x_+ (\eta) < A,
\endaligned
$$
admits an analytic extension.
We note here that we choose the branch of the square root
that is positive for $x_- < x <x_+$.

We will show that even if $R^0$ does not admit an  analytic
extension  in $\Bbb H $, the analysis of Chapter 5 in [3]
can be amended via the solution of  a scalar Riemann-Hilbert problem.

Indeed, consider the following scalar additive Riemann-Hilbert problem,
with jump on the linear segment $\Sigma = [-iA, iA]$. Let $p$ be a function
analytic in $ \Bbb C \setminus [-iA, iA]$, such that
$$
\aligned
p_+ (\eta) + p_-(\eta) =  \rho_0(\eta)= {{dR^0} \over {d\eta}}, ~~~\eta \in (-iA, iA),\\
lim_{\eta \to \infty} p(\eta) = 0.
\endaligned
$$
Here $R^0(\eta)$ is extended to the lower half of $\Sigma$
by the relation $R^0(\eta^*) = R^0(\eta)$.
The "+" side is to the left of $\Sigma$ and the
"-" side is to the right of $\Sigma$.

Note that if
$R^0$ is entire, then we can choose $p = \rho^0 = 1/2 {{dR^0} \over {d\eta}}.$
In general, our choice of initial data only ensures that $\rho^0$
is continuous.

Now, the analysis of Chapter 5 in [3] can be amended as follows.
First, let's amend the definition of $X$ in Chapter 3, which describes
the interpolant of the norming constants.
We simply set
$$
\aligned
X(\lambda) = i \pi (2K+1)  \int^{iA}_{\lambda} (p_+(\eta) +p_-(\eta) ) d\eta,
\endaligned
$$
for $\lambda $ in the linear segment $[0, iA]$.
Then, the discussion of Chapter 5 in [3], in particular from relation (5.4)
to (5.8), is  amended by substitutitng $\bar \rho^{\sigma} = p-\rho$.
More precisely, taking $\sigma = 1$,
$$
\aligned
\int_{0}^{iA}L^0_\eta(\lambda)p_-(\eta)d\eta =
\int_{C_I}
L^{C}_{\eta-}(\lambda)p(\eta)d\eta,
\endaligned
$$
and similarly, by symmetry,
$$
\aligned
\int_{-iA}^{0}L^0_\eta(\lambda)p_-(\eta^*)^*\,d\eta =
\int_{C_I^*}
L^{C}_{\eta-}(\lambda)p(\eta^*)^*\,d\eta.
\endaligned
$$
(Recall here that $L^0_\eta(\lambda)=log (\lambda-\eta),$
with a cut along the imaginary axis from $\eta$ to $-i \infty$.
In the above integral we integrate over the "-" side, while in the integral just following
we integrate over the "+" side.)
Also
$$
\aligned
\int_{0}^{iA}L^0_\eta(\lambda)p_+(\eta)d\eta =
\int_{C_F}
L^{C}_{\eta-}(\lambda)p(\eta)d\eta,
\endaligned
$$
and similarly, by symmetry,
$$
\aligned
\int_{-iA}^{0}L^0_\eta(\lambda)p_+(\eta^*)^*d\eta =
\int_{C_F^*}
L^{C}_{\eta-}(\lambda)p(\eta^*)^*d\eta.
\endaligned
$$

Next, note that
$L^{C}_{\eta+}(\lambda)=L^{C}_{\eta-}(\lambda)$ for all
$\eta\in C_I\cup C_I^*$ ``below'' $\lambda\in C_I$
and at the same time
$L^{C}_{\eta_+}(\lambda)= 2\pi i +
L^{C}_{\eta-}(\lambda)$ for $\eta\in C_I$ ``above'' $\lambda$.
This means that for $\lambda\in C$,
$$
\aligned
\int_{C} L^{C}_{\eta\pm}(\lambda)p(\eta)d\eta +
\int_{C^*} L^{C}_{\eta\pm}(\lambda)p(\eta^*)^*d\eta =
\\
\int_{C} \overline{L^{C}_\eta}(\lambda)p(\eta)d\eta
+ \int_{C^*}\overline{L^{C}_\eta}(\lambda)p(\eta^*)^* d\eta
\pm\pi i/2 \int_{C_I}p(\eta)d\eta
\pm\pi i/2 \int_{C_F}p(\eta)d\eta,
\endaligned
$$
with $\overline{L^{C}_\eta}(\lambda) = {{L^{C}_{\eta+}(\lambda)
+L^{C}_{\eta-}(\lambda) } \over 2}.$
Assembling these
results gives the expression
$$
\aligned
\tilde{\phi}(\lambda)=
\int_{C}\overline{L^{C}_\eta}(\lambda)\overline{\rho}(\eta) d\eta +
\int_{C^*}\overline{L^{C}_\eta}(\lambda)\overline{\rho}(\eta^*)^*d\eta\\
+J(2i\lambda x + 2i\lambda^2 t) -(J(2K+1)+1) ~
(\pm\pi i/2 \int_{C_I}p(\eta)d\eta
\pm\pi i/2 \int_{C_F}p(\eta)d\eta),
\endaligned
$$
valid for $\lambda \in C$, where we have introduced the
complementary density  for
$\eta\in C:
\overline{\rho}(\eta):=p(\eta)-\rho(\eta).  $
Choosing $K$ so that $J(2K+1)+1 =0$, the last term vanishes and
we simply have
$$
\aligned
\tilde{\phi}(\lambda)=\int_{C}\overline{L_\eta^{C,\sigma}}(\lambda)
\overline{\rho}(\eta)d\eta +
\int_{C^*}\overline{L_\eta^{C}}(\lambda)
\overline{\rho}(\eta^*)^*d\eta +
J(2i\lambda x + 2i\lambda^2 t).
\endaligned
$$
Compare with (5.11) of [KMM]; this  formula is less awkward,
since it does not depend on the a priori constraint
that  the contour $C$ has to go through $iA$,
a constraint that is eventually suspended anyway.

The rest of the proofs of [3] go through, with $p$ substituting $\rho^0$.
We omit the detailed discussion, but
we $do$ stress one major point on the variational problem.

As stated before in this paper,
the contour $C$  and the measure $\rho d\eta$
are characterized by a solution of a Green's  variational
problem of electrostatic kind. Indeed
$$
\aligned
E_{\phi} (\rho d\eta) = max_{C'} min_{\mu: supp(\mu)  \in C} E_{\phi} (\mu),
\endaligned
$$
where the contours $C'$   are a priori supported in the upper
half-plane minus the linear segment $[0, iA]$, and
$E_{\phi}$ is the weighted energy of a measure with respect to the external
field  given by
$$
\aligned
\phi (z) =
\int log {{ |z-\eta^*| } \over {|z-\eta|}}
\rho^0 (\eta) d\eta - Re (i\pi J \int_{z}^{iA} p (\eta) d\eta
+2iJ (z x + z^2 t) ).
\endaligned
$$
The harmonicity of $\phi$ is important to the structure of
$C, supp(\rho)$. But again, even if $\rho^0$ is not analytically extended,
it can be written as a sum of two terms that $are$.

One could write $\phi$ as
$$
\aligned
\phi (z) =
\int log {{ |z-\eta^*| } \over {|z-\eta|}}
(p_+ + p_-) (\eta) d\eta - Re (i\pi J \int_{z}^{iA} p (\eta) d\eta
+2iJ (z x + z^2 t) ).
\endaligned
$$
Again, this representation is perhaps more natural, since in setting
the variational problem it is more appropriate to think of the
"left" and "right" sides of the linear segment
$[0, iA]$ as distinct.

\bigskip

REMARK: In the main text of this paper we assumed that the solution of the variational problem does not touch the
spike $[0, iA]$ except possibly at a finite number of points. As shown in the Appendix A3, this obstacle can be
overcome by setting the variational problem on an infinite sheeted Riemann surface $\Bbb L$, where, of course,
we use the analyticity of $\rho^0$ even across the spike. Now, here we don't have that (in fact this is
the whole point of this appendix). But a careful examination of Appendix A3 shows that what we actually need is
analyticity across all but one liftings of the spike on $\Bbb L$. This we can get by simply setting our
scalar Riemann-Hilbert problem on $\Bbb L$ and letting the jump be a single copy of the spike $[0, iA]$
in $\Bbb L$. The scalar Riemann-Hilbert problem on $\Bbb L$ can be explicitly solved by mapping conformally
$\Bbb L$ to $\Bbb C$.

\bigskip

CONCLUSION:
The moral of the story is that if $\rho^0$ does not admit a holomorphic  extension,
we can write it as the average of two functions $p_-, p_+$ that can be extended
to the left and right of the segment $[0, iA]$ respectively,
and proceed as before, with $\rho^0$  substituted by $p$.

\newpage

7. REFERENCES

\bigskip

[1]  A. A. Gonchar and E. A. Rakhmanov, Equilibrium Distributions and Degree
of Rational Approximation of Analytic Functions,
Math. USSR Sbornik, v. 62, pp.305--348, 1989.

[2] E. B. Saff, V. Totik, Logarithmic Potentials
with External Fields, Springer Verlag, 1997.

[3] S. Kamvissis, K. T.-R. McLaughlin, P. D. Miller, Semiclassical Soliton
Ensembles
for the Focusing Nonlinear Schr\"odinger Equation,
Annals of Mathematics Studies, v.154, Princeton University Press, 2003.

[4] P. Deift, X.Zhou, A Steepest Descent Method for
Oscillatory Riemann-Hilbert Problems, Annals of Mathematics, v.137, n.2, 1993,
pp.295-368.

[5] E. A. Perevozhnikova and E. A. Rakhmanov, Variations of the Equilibrium
Energy and S-property of Compacta of Minimal Capacity, preprint, 1994.

[6] J. Dieudonn\'e, Foundations of Modern Analysis, Academic Press, 1969.

[7] N. S. Landkof, Foundations of Modern Potential Theory, Springer Verlag,
1972.

[8] G. M. Goluzin, Geometric Theory of Functions of a Complex Variable,
Translations of Mathematical Monographs, v.26, AMS 1969.

\end{document}